\documentclass{article}
\usepackage{amssymb}

\begin{document}

\vskip2cm
\newcommand{\gothS}{\digamma}
\newcommand{\gothD}{\mho}
\newcommand{\rl}{{\bf R}^1}
\newcommand{\rr}{{\bf R}}
\newcommand{\rk}{{\bf R}^k}
\newcommand{\rn}{{\bf R}^n}
\newcommand{\ri}{{\bf R}^i}
\newcommand{\nn}{{\bf N}}
\newcommand{\zplus}{{\bf Z}_{+}}
\newcommand{\zz}{{\bf Z}}
\newcommand{\rntn}{{\bf R}^{n \times n}}
\newcommand{\rnextnex}{{\bf R}^{n \times n}}
\newcommand{\rnextpex}{{\bf R}^{n \times p}}
\newcommand{\rnex}{{\bf R}^{n }}
\newcommand{\rpex}{{\bf R}^{p}}
\newcommand{\rktk}{{\bf R}^{k \times k}}
\newcommand{\Linfty}{L_{\infty}(I;\rl)}
\newcommand{\Ll}{L_1(I;\rl)}
\newcommand{\Cheb}{C(I;\rl)}
\newcommand{\Chebk}{C(I;\rk)}
\newcommand{\Gchebk}{C^1(I;\rk)}
\newcommand{\Gcheb}{C^1(I;\rl)}
\newcommand{\zt}{\zeta}
\newcommand{\lmb}{({\lambda}_1,...,{\lambda}_k{)}^T}
\newcommand{\nrho}{\rho}
\newcommand{\neta}{\eta}
\newcommand{\rbn}{{\bf R}^{N}}
\newcommand{\xov}{{\tilde x}^{0}}
\newcommand{\tauex}{{\overline \tau }}
\newcommand{\ff}{{F }}
\newcommand{\rtwo}{{\bf R}^2}
\newcommand{\rthree}{{\bf R}^3}
\newcommand{\rfour}{{\bf R}^4}
\newcommand{\rmnu}{{\bf R}^{m_{\nu}}}
\newcommand{\rntm}{{\bf R}^{n \times m}}
\newcommand{\rmiplusl}{{\bf R}^{m_{i+1}}}
\newcommand{\rmlplusmi}{{\bf R}^{m_1+ \ldots +m_i}}
\newcommand{\rmnuplusl}{{\bf R}^{m_{\nu + 1}}}
\newcommand{\rkplusl}{{\bf R}^{k + 1}}
\newcommand{\rmlplusdp}{{\bf R}^{m_1+ \ldots +m_p}}
\newcommand{\rmpplust}{{\bf R}^{m_{p+2}}}

\newcommand{\rml}{{\bf R}^{m_1}}
\newcommand{\rmt}{{\bf R}^{m_2}}
\newcommand{\rmpmin}{{\bf R}^{m_{p-1}}}
\newcommand{\rsumpmin}{{\bf R}^{m_1 + ...+m_{p-1}}}
\newcommand{\rmtm}{{\bf R}^{m \times m}}
\newcommand{\rmtrmpl}{{\bf R}^{m \times m_{p+1}}}
\newcommand{\rmp}{{\bf R}^{m_p}}
\newcommand{\rmm}{{\bf R}^{m}}
\newcommand{\rmpplusl}{ {\bf R}^{m_{p+1}} }
\newcommand{\rmppluslzpt}{ {\bf R}^{m_{p+1}} }
\newcommand{\rmpplusldt}{ {\bf R}^{m_{p+1}} }
\newcommand{\rmitmj}{ {\bf R}^{m_{i} \times m_{j}} }
\newcommand{\rktrmpl}{{\bf R}^{k \times m_{p+1}}}
\newcommand{\rktkplusl}{{\bf R}^{k \times (k+1)}}
\newcommand{\rmltmtml}{ {\bf R}^{m_{1} \times {(m_{2}-m_1)} }   }
\newcommand{\rmltml}{ {\bf R}^{m_{1} \times m_{1} }   }
\newcommand{\rmtminml}{ {\bf R}^{m_{2}-m_{1}} }
\newcommand {\rmminmp}{ {\bf R}^{m-m_{p}} }
\newcommand {\rrmmpp}{ {\bf R}^{ ({m-m_{p}})\times m_{p} }}
\newcommand{\rmtminm}{{\bf R}^{m_{2} -m} }
\newcommand {\rkminmp}{ {\bf R}^{k-m_{p}} }
\newcommand {\rrkmpp}{ {\bf R}^{ ({k-m_{p}})\times m_{p} }}
\newcommand{\rmtmink}{{\bf R}^{m_{2} -k} }
\newcommand{\rmi}{{\bf R}^{m_i}}
\newcommand{\rmj}{{\bf R}^{m_j}}
\newcommand{\rktm}{{\bf R}^{k \times m}}
\newcommand{\nxi}{\xi}
\newcommand{\rmlplusmp}{{\bf R}^{m_1 + ...+m_{p}}}

\newcommand{\ekk}{E^{2k+2}}

\newcommand{\mll}{\parallel   w   (\cdot)   {\parallel}_{C ([t_0,T];\rpex) }}
\newcommand{\ml}{{\parallel   w   (\cdot)   {\parallel}_{C ([t_0,T];\rpex) }^{-1}}}
\newcommand{\dlm}{{\delta \lambda} = ( {\delta \lambda}_1,\ldots,{\delta\lambda}_k) \in { \bf R}^k}
\newcommand{\wdlms}{w(s)}
\newcommand{\wdlmt}{w(\tau)}

\newcommand{\uhlms}{{\hat u}(s)+{\sum_{j=1}^{k}{ {\lambda_j}{\hat w}_j(s)}}}
\newcommand{\uhlmt}{{\hat u}(\tau)+{\sum_{j=1}^{k}{ {\lambda_j}{\hat w}_j(\tau)}}}

\newcommand{\uhtlms}{{\hat u}_{\lambda}(\cdot)}
\newcommand{\uhtlmt}{{\hat u}_{\lambda}(\tau)}

\newcommand{\wblms}{{\sum_{j=1}^{k}{ {\bar {\lambda_j}}{w}_j(s)}}}
\newcommand{\wblmt}{{\sum_{j=1}^{k}{ {\bar {\lambda_j}}{w}_j(\tau)}}}

\newcommand{\whblms}{{\sum_{j=1}^{k}{ {\bar {\lambda_j}}{\hat {w}_j}(s)}}}
\newcommand{\whblmt}{{\sum_{j=1}^{k}{ {\bar {\lambda_j}}{\hat {w}_j}(\tau)}}}

\newcommand{\wmbblmt}{\hat {w_{\bar{\lambda}}}(\cdot)}
\newcommand{\whbblms}{\hat {w_{\bar{\lambda}}}(s)}
\newcommand{\whbblmt}{\hat {w_{\bar{\lambda}}}(\tau)}
\newcommand{\whbblm}{\hat {w_{\bar{\lambda}}}}

\newcommand{\dfdx}{\frac {\partial f}{\partial x}}
\newcommand{\dfidx}{\frac {\partial {f_i}}{\partial x}}
\newcommand{\dgdx}{\frac {\partial g}{\partial x}}
\newcommand{\dgidx}{\frac {\partial {g_i}}{\partial x}}

\newcommand{\dfdu}{\frac {\partial f}{\partial u}}
\newcommand{\dfidu}{\frac {\partial {f_i}}{\partial u}}
\newcommand{\dgdu}{\frac {\partial g}{\partial u}}
\newcommand{\dgidu}{\frac {\partial {g_i}}{\partial u}}

\newcommand{\ulms}{\overline v (s) }
\newcommand{\ulmt}{\overline v (\tau)  }
\newcommand{\udlms}{{\overline v} (s) + w (s)}
\newcommand{\udlmt}{\overline v (\tau ) + w (\tau)}
\newcommand{\utdlms}{\overline v (s) + {\bar {\theta}_i} (s)  w(s)}
\newcommand{\utdlmt}{\overline v (\tau) + {\bar   {\vartheta}_i }  (s,\tau) w(\tau)}

\newcommand{\xlms}{x_{\overline v (\cdot )}(x^0,  s)}
\newcommand{\xlmt}{x_{\overline v (\cdot )}(x^0,   {\tau})}

\newcommand{\xhlm}{\hat {x_{\lambda}}}

\newcommand{\xhblms}{\hat {x_{\bar{\lambda}}}(s)}
\newcommand{\xhblmt}{\hat {x_{\bar{\lambda}}}({\tau})}
\newcommand{\xhblm}{\hat {x_{\bar{\lambda}}}}
\newcommand{\xdlms}{x_{\overline v (\cdot ) + w(\cdot)}( x^0 ,s)}
\newcommand{\xdlmt}{x_{\overline v (\cdot ) + w(\cdot)}( x^0 ,\tau)}
\newcommand{\xlm}{x_{\overline v (\cdot)} }
\newcommand{\xdlm}{x_{ \overline v (\cdot) + w(\cdot) } }
\newcommand{\xilm}{x_{\overline v (\cdot) }^i}
\newcommand{\xidlm}{x_{\overline v (\cdot)  + w(\cdot) }^i}
\newcommand{\xtdlms}{x_{\overline v (\cdot) } (x^0,s) + {{\theta}_i(s)} (x_{\overline v (\cdot) + w(\cdot) }(x^0,s) - x_{\overline v (\cdot)}(x^0,s))}
\newcommand{\xtdlmt}{x_{\overline v (\cdot)}(x^0,\tau) + {{\vartheta}_i(s,\tau )}(x_{\overline v (\cdot) + w (\cdot)}(x^0,\tau) - x_{\overline v (\cdot)}(x^0,\tau))}
\newcommand{\zidlmlm}{z_{{w(\cdot)},{\overline v (\cdot)}}^i}
\newcommand{\zdlmlm}{z_{{w(\cdot)},{\overline v (\cdot)}}}
\newcommand{\zbdlmlm}{z_{{\bar {\lambda}},{\lambda +{\delta \lambda}}}}
\newcommand{\zblm}{z_{{\bar {\lambda}},{\lambda}}}


\begin{large}

\newcommand{\n}[1]{\refstepcounter{equation}\label{#1}
\eqno{(\arabic{equation})}}
\renewcommand{\theequation}{\arabic{equation}}

\begin{Large}
\begin{center}
{\bf The solution of the global controllability problem for
the triangular systems in the singular case}
\end{center}

$$ \; $$

\begin{center}
Valery I. Korobov, and  Svyatoslav S. Pavlichkov
\end{center}

\end{Large}
$$\; $$

$$ \; $$
\begin{large}
\begin{center}
{\bf Abstract.}
\end{center}
The solution of the global controllability problem is
obtained for a class of the triangular systems of O.D.E.
that are not feedback linearizable. The introduced class is
a generalization of the classes of triangular systems
investigated before. The solution of the problem is based
on the approach proposed in \cite{ssp_op1} for the
triangular systems of the Volterra equations. This yields
the same properties of the considered class of triangular
systems as those established in \cite{ssp_op1} for the
Volterra systems. As well as in \cite{ssp_op1}, for the
current class of triangular systems, it is proven that
there exists a family of continuous controls that
   solve the global controllability problem for the considered class
 and continuously depend on the initial and the terminal states.
 As well as in \cite{ssp_op1}, this implies the
 global controllability of the bounded perturbations of
 the current class. In contrast with \cite{ssp_op1},
 to prove the existence of the desired family of open-loop
 controls, we construct a family of closed-loop ones each
 of which steers the corresponding initial state into
 an appropriate neighborhood  of an appropriate terminal
 point.

\end{large}
$\qquad $

{\bf Key words:} Nonlinear control, triangular form, global
controllability, feedback linearization.

{\bf Mathematics subject classification:} 93C10, 93B10,
93B11, 93B05, 93B52.

$\qquad$

\begin{center}
{\bf 1. Introduction and the statement of the main results.
}
\end{center}

In this paper, we consider a control system $$ \dot x(t) =
f(t,x(t),u(t)), \; \; \; \; \; \; \; \; t \in [t_0,T],
\n{r4_1}$$
 where $x {\in} \rn$ is the state,
$u {\in} \rmm$ is the control, and function $f$ has the
following "triangular" form $$f(t,x,u) = \left(
\begin{array}{l} f_1(t,x_1,x_2) \\ f_2(t,x_1,x_2,x_3) \\
 \quad \ldots \ldots \ldots\quad \\
f_{\nu-1} (t,x_1,x_2,...,x_{\nu})\\ f_{\nu}
(t,x_1,x_2,...,x_{\nu},u)
\end{array}\right),\; \n{r4_2}$$
 and
$$ \nu \in \nn;\; \; \;   x_i \in \rmi, \; \; \;
f_i(t,x_1,...,x_{i+1}) \in \rmi, \; \; \; i = 1,...,\nu;\;
\;u \in \rmm=\rmnuplusl; \; \; \;  m{=}m_{\nu{+}1};$$ $$
m_i {\leq} m_{i{+}1}, \; \; i{=}1,...,\nu;    \; \; \;   n
{=} m_1 {+} ...{+} m_{\nu};\; \; \;
x{=}(x_1,...,x_{\nu}{)}^T \in \rn {=} \rml
{\times}\rmt{\times}{\ldots}{\times}\rmnu.$$ We assume that
function $f$ satisfies the following conditions:

(I)  $f,$ $\frac{\partial f}{\partial x},$ and
 $\frac{\partial f}{\partial u}$ are of classes
 $C([t_0,T]\times\rn \times\rmm; \rn),$
  $ C([t_0,T] \times \rn \times\rmm; \rntn),$
{\sl and}  $ C([t_0,T] \times \rn \times \rmm; \rntm)$ {\sl
respectively.}

    (II) {\sl For each $ i{=} 1,...,\nu, $ and  each
 $(t,x_1,...,x_i) {\in} [t_0,T] {\times} \rml {\times} ...{\times} \rmi,$ the map
 $f_i (t,x_1,...,x_i,\cdot)$ is of class
 $C^{m_{i+1} - m_i +1} (\rmiplusl;\rmi),$ and
 $ f_i (t,x_1,...,x_i,\rmiplusl)= \rmi.$}


Triangular systems appeared at the earliest stage in the
development of the nonlinear control theory
 -- \cite{kor1}.
Triangular form arose initially from a specific problem
concerned with control of satellites (\cite{kovalev2}),
but in the end it became  fruitful in modeling various
physical and engineering systems (see, for instance
\cite{borisov}, \cite{kovalev3}, or related works
\cite{fliess}, \cite{murray}).

One of the causes of this was the habitual situation when
the output of a given control system affects the input of
another control system, which produces
 the triangular form - \cite{borisov},
\cite{kovalev3}. This fact  eventually led to the general
concept of flatness \cite{fliess} -- the last being a very
popular generalization of the notion of exact
linearization. On the other hand, in the case of the
cascade form, the so-called backstepping technique is
employed for asymptotic stabilization very often --
\cite{Fre_Koko}, \cite{kokotovic_sussmann},
\cite{lin_kanellakopoulos}, \cite{saberi}. In addition,
triangular systems appear naturally in the general exact
linearization approach (see \cite{Cel_Nij},
\cite{jakubczyk}, \cite{kor1}, \cite{Resp}), which has a
large number of applications -- \cite{atanasov},
\cite{D_Andrea}, \cite{fliess}, \cite{kovalev3},
\cite{MRicc}, \cite{kovalev2}, \cite{singh1}.

In most works cited above, the triangular form is treated
as a mere specific form of the dynamics of a feedback
linearizable system however. (The exception is the papers
where the main objective is only the design of stabilizers
without finding the appropriate states coordinates --
\cite{Fre_Koko}, \cite{kokotovic_sussmann}, \cite{saberi},
\cite{tsinias}). More precisely, the standard requirement
for the triangular form (\ref{r4_2}) is that the conditions
$|\frac{\partial f_i}{\partial x_{i+1}}| \not=0,$
$i=1,...,\nu$ $(x_i {\in} \rl),$ or, even, $|\frac{\partial
f_i}{\partial x_{i+1}}| {\geq} a{>}0,$ $i=1,...,\nu,$ hold.
This condition implies the feedback equivalence of the
system to the linear canonical form $\dot z_i = z_{i+1},$
$i{=}1,...,\nu{-}1,$ $\dot z_{\nu} =v;$ and, conversely, if
this condition does not hold, then the system may be not
feedback linearizable in the corresponding domain (see
example 1.1 below).

The case when these conditions do not hold, which is called
the "singular" case, is considered in works \cite{Cel_Nij},
\cite{Resp}. In these papers, the triangular systems are
studied provided that the set of the regular points (i.e.,
the set of all $x$ such that $|\frac{\partial f_i
(x)}{\partial x_{i+1}}| \not=0,$ $i{=}1,...,\nu{-}1,$) is
open and dense in the state space however. This may be not
the case in some simple examples (see again example 1.1
which is given below). That is why,
 as it is concluded in \cite{Cel_Nij},
the singular case requires further investigation.

On the other hand, most results concerned with the concept
of feedback linearization are essentially local (again, as
the exception, we can mention work \cite{Resp} which we are
going to generalize in some ways). This is another motive
for introducing new classes of the triangular systems that
are wider than those investigated before \cite{kor1},
\cite{pavl1}, \cite{Resp}, \cite{Cel_Nij}.

Therefore, our main goal is to find a generalization, of
the triangular form, that can be studied globally. The
class of triangular systems that is defined by our
conditions (I),(II) is such a generalization. In this
paper, we restrict our study to the solution of the global
controllability problem  for this class. For this, we use
the approach that was developed in \cite{pavl3},
\cite{ssp_op1} for the triangular systems of the Volterra
equations. As the outcome, we obtain that the formulation
of the main theorems are the same both in the current work
and in \cite{ssp_op1} (but, of course, the classes of the
control systems that are considered in \cite{ssp_op1} and
in the current work differ essentially).

Like in \cite{ssp_op1}, for (\ref{r4_1}), we construct a
family of open-loop continuous controls parametrized by the
initial and the terminal states such that every element of
this family steers the corresponding initial state to the
corresponding terminal one and continuously depends on the
initial and the terminal states (theorem 1.1.). This
automatically yields some kind of robustness; in
particular, we obtain the global controllability of the
bounded perturbations of system (\ref{r4_1}) (theorem 1.2)


    Following \cite{ssp_op1}, we consider a perturbation of system
 (\ref{r4_1}) of the  form
$$ \dot x(t) = f(t,x(t),u(t)) + h(t,x(t),u(t)),\; \; \; \;
\;  t \in [t_0,T], \n{r4_3}$$ where $h$ satisfies the
conditions:

 (III).  $h {\in} C([t_0,T]{\times}\rn{\times}\rmm; \rn),$
{\sl and for each compact set }
  $Q {\subset} \rn {\times} \rmm,$ {\sl there exists} $L_Q {>}0$ {\sl such that, for each}
 $t {\in} [t_0,T],$ {\sl each }
 $(x^1,u^1) {\in} Q,$ {\sl and each }
 $(x^2,u^2) {\in} Q,$ {\sl we obtain:}
 $$ |h(t,x^1,u^1) - h(t,x^2,u^2) | \leq L_Q (|x^1 - x^2| +|u^1-u^2|). $$

  (IV).  {\sl There exists  $H{>}0$  such that, for each
  $(t,x,u) {\in} [t_0,T]{\times}\rn{\times}\rmm,$ we have:
$|h(t,x,u)| {\leq} H.$}

 Throughout the paper, for each  $x^0 {\in} \rn,$ each  $ u(\cdot) {\in} L_{\infty}([t_0,T]; \rmm),$
and each  $\tau{\in}[t_0,T],$ by
 $t \mapsto x(t,\tau, x^0, u(\cdot))$ we denote
the trajectory, of system
 (\ref{r4_1}),  that is defined by the control
 $u(\cdot)$ and by the initial condition
 $x(\tau){=}x^0$ on some maximal subinterval
 $ J_1{\subset}[t_0,T],$ $\tau{\in}J_1.$ The main results of the current
work are the following theorems 1.1-1.3.

 {\bf  Theorem 1.1.} {\sl  Assume that, for system (\ref{r4_1}), function
$f$ has triangular form
 (\ref{r4_2}) and satisfies conditions
 (I) and  (II).
Then, there exists a family of controls
   $\{ u_{(x^0,x^T)} (\cdot) {\}}_{(x^0,x^T) \in \rn \times \rn}$
such that the map  given by $(x^0,x^T) \mapsto
u_{(x^0,x^T)} (\cdot)$ is of class $C(\rn \times\rn;
C([t_0,T]; \rmm)),$ and, for every
 $(x^0,x^T) \in \rn \times \rn,$    the trajectory
 $t {\mapsto} x(t,t_0,x^0,u_{(x^0,x^T)} (\cdot))$ is defined
 for all $t{\in}[t_0,T],$  and
 $x(T,t_0,x^0,u_{(x^0,x^T)} (\cdot)) = x^T.$
  }

  {\bf   Theorem 1.2.}   {\sl  Assume that, for system (\ref{r4_1}),
 function   $f$ has  triangular form
  (\ref{r4_2}), satisfies
   (I), (II),  and satisfies the global Lipschitz condition w.r.t. $x$
and $u,$ i.e., there exists
   $L>0$   such that, for each   $t{\in}[t_0,T],$
each $(x^1,u^1) {\in} \rn {\times} \rmm,$ and each
   $(x^2,u^2) \in \rn \times \rmm,$ we have
   $$ |f(t,x^1,u^1) - f(t,x^2,u^2)| \leq L(|x^1-x^2| + |u^1 - u^2|).  $$
Suppose that function   $h$ satisfies (III) and (IV). Then
system    (\ref{r4_3}) is globally controllable in time
    $[t_0,T]$    by means of controls from class  $C([t_0,T]; \rmm).$
    }

  {\bf  Theorem 1.3.}   {\sl Assume that  $f$     has triangular form
    (\ref{r4_2})
and satisfies   (I), (II).    Then system  (\ref{r4_1})
  is globally controllable in time $[t_0,T]$ by means of controls from class
  $C([t_0,T]; \rmm).$
    }

  {\bf  Example 1.1.} Consider the system
    $$ \left\{ \begin{array}{l}
 \dot x_1 (t) = g(x_2(t)) \\
 \dot x_2 (t) = u(t)  \end{array}\right.   \; \; \;\; \; \;  \; t \in [0,T], \n{pr1} $$
where $(x_1,x_2)^T {\in} \rl {\times} \rl {=} \rtwo$ is the
state,  $u {\in} \rl$ is the control, $T>0$ is an arbitrary
fixed number, and $g(\cdot)$ is given by $$ g(y) = \left\{
\begin{array}{l}
 0  \; \; \;\; \; \; \; \; \; \; \; \; \; \; \; \;\; \; \; \; \;   \; \; \;\;\;\;\; \;  \mbox{ if } y \leq 2 \\
  (y-2{)}^2 \sin (y-2)  \; \; \; \mbox{ if } y > 2 \end{array}\right.    \n{pr10} $$
It is clear that system (\ref{pr1}) has triangular form
(\ref{r4_2}) and satisfies conditions  (I) and (II).

  Let us prove that (\ref{pr1})
is not globally feedback equivalent to the canonical
system, i.e.,  there are no functions $\phi_1(t,x_1,x_2),$
$\phi_2(t,x_1,x_2),$ and $\phi_3(t,x_1,x_2,u),$ of class
$C^1,$ such that (A) for each
   $ t {\in} [0,T],$ the maps
$(x_1,x_2{)}^T {\mapsto} (\phi_1(t,x_1,x_2),
\phi_2(t,x_1,x_2) {)}^T$   and $(x_1,x_2,u{)}^T \mapsto
(\phi_1(t,x_1,x_2),\phi_2(t,x_1,x_2),
\phi_3(t,x_1,x_2,u){)}^T$ are diffeomorphisms of $\rtwo$
onto $\rtwo,$ and of $\rthree$ onto $\rthree$ respectively,
and (B) if
  $(x_1(\cdot), x_2(\cdot){)}^T $ is a trajectory of
 (\ref{pr1})
with some control  $u(\cdot),$ then the functions
$v(\cdot)$  and  $(z_1(\cdot),z_2(\cdot){)}^T $ that are
defined by $$ \begin{array}{l}
  z_1 (t) = {\phi}_1(t,x_1(t),x_2(t)) \\
 z_2 (t) = {\phi}_2(t,x_1(t),x_2(t))\\
v(t) = {\phi}_3(t,x_1(t),x_2(t),u(t)),  \end{array}   \; \;
\;\; \; \;  \; t \in [0,T],  \n{pr4}$$ satisfy $$\left\{
\begin{array}{l}
 \dot z_1 (t) = z_2(t) \\
 \dot z_2 (t) = v(t)  \end{array}\right.   \; \; \;\; \; \;  \; t \in [0,T], \n{pr5} $$
and conversely, if  $(z_1(\cdot),z_2(\cdot){)}^T$ is a
trajectory of system (\ref{pr5}) with some control
$v(\cdot),$ then the functions $u(\cdot)$ and
$(x_1(\cdot),x_2(\cdot){)}^T$ that are defined by
(\ref{pr4}) satisfy (\ref{pr1}).

Assume the converse, i.e., that such a global feedback
transformation does exist. Put: $x_1^{\ast} {=}x_2^{\ast}
{=} u^{\ast} {=} 0.$
 Let $(x_1(\cdot),x_2(\cdot){)}^T$ be the trajectory, of
 (\ref{pr1}), that is defined by the control
 $u(t) {=} u^{\ast},$   $t{\in}[0,T],$
and by the initial condition
  $(x_1(0),x_2(0){)}^T {=} (x_1^{\ast},x_2^{\ast}{)}^T.$
Then  $(x_1(t),x_2(t){)}^T {=}( x_1^{\ast},
x_2^{\ast}{)}^{T} {=} 0 {\in} \rtwo$
 for all  $t {\in} [0,T].$ By $(z_1(\cdot),z_2(\cdot){)}^T$
and  $v(\cdot)$  denote the trajectory and the control, of
system (\ref{pr5}), that are defined by (\ref{pr4}). Put:
$z_1^{\ast}= \phi_1(0,x_1^{\ast},x_2^{\ast}),$ $z_2^{\ast}
= \phi_2(0,x_1^{\ast},x_2^{\ast}).$ Then, it is clear that
the controls
 $v_1(t) {:=} \frac{6}{T^2} {-}\frac{12 t}{T^3}$
and  $v_2(t) {=} - \frac{2}{T} {+} \frac{6t}{T^2},$ $t
{\in} [0,T],$ steer the origin $(0,0{)}^T {\in} \rtwo$ into
vectors $e_1 {=} (1,0{)}^T$ and $e_2 {=} (0,1{)}^T$
respectively in time $[0,T]$ w.r.t.  (\ref{pr5}). For each
$\mu {=} ({\mu}_1,{\mu}_2{)}^T {\in} \rtwo,$ by $v_{
({\mu}_1,{\mu}_2) } (\cdot) $
 denote the control  $v_{({\mu}_1,{\mu}_2)} (\cdot) {=} v(\cdot) {+} {\mu}_1 v_1(\cdot) {+}{\mu}_2 v_2(\cdot),$
and let $(z_1({\mu}_1, {\mu}_2, \cdot),
z_2({\mu}_1,{\mu}_2,\cdot){)}^T$ be the trajectory, of
(\ref{pr5}), defined by the control $v_{(\mu_1,\mu_2)}
(\cdot)$ and by the initial condition
$(z_1(\mu_1,\mu_2,0),z_2(\mu_1,\mu_2,0){)}{=} (z_1^{\ast},
z_2^{\ast}).$
 Let  $(x_1(\mu_1, \mu_2, \cdot), x_2(\mu_1,\mu_2,\cdot){)}^T$
and $u_{(\mu_1,\mu_2)} (\cdot)$ be defined by $$z_i (\mu_1,
\mu_2,t )=\phi_i (t, x_1(\mu_1,\mu_2,t),
x_2(\mu_1,\mu_2,t)), \; \; i = 1,2; $$ $$v_{(\mu_1,\mu_2)}
(t) = \phi_3 (t,x_1(\mu_1,\mu_2,t), x_2(\mu_1,\mu_2,t),
u_{(\mu_1,\mu_2)} (t));\; \;  t \in [0,T]. \n{pr6}$$ Then,
since  (\ref{pr5}) is a linear system, we obtain that
$$\left(z_1(\mu_1,\mu_2,T),z_2(\mu_1,\mu_2,T)\right) =
\left( z_1(T) +\mu_1, z_2(T) +\mu_2 \right),$$ i.e., the
map $(\mu_1,\mu_2{)}^T {\mapsto}
(z_1(\mu_1,\mu_2,T),z_2(\mu_1,\mu_2,T){)}^T$ is a
diffeomorphism of $\rtwo$ onto $\rtwo;$ therefore, since
(\ref{pr4}) satisfies condition (A),
 from
 (\ref{pr6}), it follows that
$(\mu_1,\mu_2{)}^T \mapsto
(x_1(\mu_1,\mu_2,T),x_2(\mu_1,\mu_2,T){)}^T$ is (in
particular) a diffeomorphism of some neighborhood of
 $(\mu_1,\mu_2{)}^T {=} 0 {\in} \rtwo$
onto some neighborhood of $(x_1^{\ast},x_2^{\ast}) {=} 0
{\in} \rtwo.$

On the other hand, from (\ref{pr6}), it follows that the
maps $(\mu_1,\mu_2{)}^T \mapsto u_{(\mu_1,\mu_2)} (\cdot)$
and $(\mu_1,\mu_2{)}^T \mapsto
(x_1(\mu_1,\mu_2,\cdot),x_2(\mu_1,\mu_2,\cdot){)}^T$ are of
classes $C(\rtwo; C([0,T];\rl))$ and $C(\rtwo;
C([0,T];\rtwo))$ respectively; hence, taking into account
that $$(x_1(0,0,t),x_2(0,0,t){)}^T {=}
(x_1^{\ast},x_2^{\ast}{)}^T {=} 0 {\in} \rtwo, \;  \; \;
\;\quad \; u_{(0,0)} (t) {=} u^{\ast} {=} 0,\; \; \quad\;
\;  \; \; \; t {\in} [0,T],$$ we get the existence of a
neighborhood $B$ of point $0 {\in} \rtwo$ such that, for
each $(\mu_1,\mu_2{)}^T {\in} B,$ we obtain $$-2 <
x_2(\mu_1,\mu_2,t) <2 \; \; \; \; \;   \;  \mbox{ for all }
t {\in} [0,T]. \n{pr7}$$ By the construction,
$(x_1(\mu_1,\mu_2,\cdot),x_2(\mu_1,\mu_2,\cdot){)}^T$ is
the trajectory, of (\ref{pr1}), defined by the control
$u_{(\mu_1,\mu_2)} (\cdot)$ and by the initial condition
    $x_1(\mu_1,\mu_2,0) {=} x_2(\mu_1,\mu_2,0) {=} 0;$
therefore, from (\ref{pr7}),  and from (\ref{pr10}), we
get: $x_1 (\mu_1,\mu_2, t) = 0 $ whenever $t {\in} [0,T]$
and $(\mu_1,\mu_2{)}^T {\in} B.$ Hence, the rank of the map
$(\mu_1,\mu_2{)}^T {\mapsto} (x_1(\mu_1,\mu_2,T),
x_2(\mu_1,\mu_2,T){)}^T$ is less than 2 at each
$({\mu}_1,{\mu}_2) \in B.$ This contradicts the fact that
this map is a local diffeomorphism at $\mu = 0$ and proves
that system (\ref{pr1}) is not globally feedback equivalent
to system (\ref{pr5}).

Let us remark that, system (\ref{pr1}) satisfies conditions
(I), (II) so that we can apply theorems 1.1, 1.3 to
(\ref{pr1}). On the other hand, since $g(\cdot)$ is given
by (\ref{pr10}), system (\ref{pr1}) does not satisfy the
global Lipschitz condition, and, therefore, we can not
apply theorem 1.2 to the bounded perturbations of
(\ref{pr1}). However, if $g(\cdot)$   were given by $$ g(y)
= \left\{ \begin{array}{l}
 0  \; \; \;\; \; \; \; \;  \; \; \;\; \; \; \; \;  \; \; \;\; \; \; \; \;  \; \; \; \; \; \;\; \; \; \; \;   \; \; \;\;\;\;\; \;  \mbox{ if } y \leq 2 \\
 {\rm ln}^2 (y-1) \times \sin\{{\rm ln} (y-1) \} \; \; \; \mbox{ if } y > 2, \end{array}\right.     $$
then we could apply all our theorems 1.1-1.3 to system
(\ref{pr1}), but (\ref{pr1}) would not be feedback
equivalent to
 (\ref{pr5}) as well.

Before proving theorems 1.1-1.3, we must stress the
following three main points.

{\bf Remark 1.1.} The proof of theorem 1.2 follows from
theorem 1.1 and from the Brouwer fixed-point theorem. We
omitt the argument, which is the same as  in
\cite{ssp_op1}, (see sect.2). It is clear that theorem 1.3
is a corollary of theorem 1.1 as well.  Therefore, the
proof of theorem 1.1 is the only goal of the next sections
of the current issue.

{\bf Remark 1.2.} As we mention above, the proof of theorem
1.1 is based on the same approach as in \cite{ssp_op1}.
However, one encounters the following two problems on this
way.

A) Condition (ii) from \cite{ssp_op1} (see sect.1) implies
the complete controllability of the linearized control
system (system (24) in \cite{ssp_op1}) around every
trajectory of the triangular system. In other words, for
the class of systems considered in \cite{ssp_op1}, the
input-output map $u(\cdot) \mapsto x(T,t_0,x^0,u(\cdot))$
is of the full rank at $u(\cdot)$ whatever $x^0$ and
$u(\cdot)$ are chosen. Condition (II) from the current work
does not ensure the same property for our system
(\ref{r4_1}). For instance, if a trajectory of system
(\ref{pr1}) from our example 1.1 lies in the set
$\{(x_1,x_2{)}^T \in \rtwo |\; \;  x_2 \leq 2\},$ then the
linearized control system around such a trajectory will not
be completely controllable.

B)    For the triangular systems from \cite{ssp_op1},
conditions (ii) and  (iii) guarantee much more than only
the existence of the solution for the last integral
equation of the triangular system w.r.t. the control (see
Eq.(8) from \cite{ssp_op1}); the obtained controls
continuously depend on the corresponding parameter $\chi$
(w.r.t the metric of $C([t_0,T]; \rl)$). And again,
condition (II) from the current work does not imply any
similar property for (\ref{r4_1}). Assume, for simplicity,
that $m_i {=}1,$ $i{=}1,...,\nu;$ fix some $i{=}1,...,\nu,$
and let $(\zeta,\xi){\mapsto} x_j (\zeta,\xi, \cdot),$
$j{=}1,...,i,$ be arbitrary maps of class $C(\ri \times
\ri; C^1([t_0,T];\rl)).$ Given $(\zeta,\xi) {\in}
\ri{\times}\ri,$ find a measurable function
$v_{(\zeta,\xi)} (\cdot)$ such that $$\dot x_i
(\zeta,\xi,t) = f_i
(t,x_1(\zeta,\xi,t),...,x_i(\zeta,\xi,t), v_{(\zeta,\xi)}
(t)) \; \; \; \; \; \mbox{ whenever } \; \; t
{\in}[t_0,T].$$ Although the existence of such
$v_{(\zeta,\xi)} (\cdot)$ can be proved via a slight
modification of the well-known Filippov lemma (see
\cite{pavl1}) the obtained map $(\zeta,\xi) {\mapsto}
v_{(\zeta,\xi)} (\cdot)$ may  be discontinuous even as a
map to
 $L_1([t_0,T]; \rmiplusl).$

{\bf Remark 1.3.} To cope with the above-mentioned problem
B), we prove lemma 3.4 (see sect.3), which is the main
difference between the current proof and that of
\cite{ssp_op1}. We emphasize that, although our original
problem is to obtain appropriate open-loop controls, we
construct a feedback law while proving lemma 3.4. - see
(\ref{r42_55_a}), and (\ref{r42_74}). This construction not
only yields the desired robustness properties of the
obtained controls but also is a method for solving the
problem of global stabilization for our class of systems.


\vskip10mm
\begin{center}
{\bf 2. The reduction of the main result to the problem of
the controllability with boundary conditions imposed on
controls.}
\end{center}

Choose  arbitrary $t_1 {\in} ]t_0,T[,$ and  $x_1^{\ast}
{\in} \rml.$ From (II), and from the Sard theorem we get
(by induction on $i{=}1,...,\nu$)
 the existence of $z_i^{\ast}{\in}\rmi$ and
 $x_{i{+}1}^{\ast}{\in}\rmiplusl,$
$i{=}1,...,{\nu},$ such that $$ z_i^{\ast} =
f_i(t_1,x_1^{\ast},..., x_{i+1}^{\ast}),  \; \mbox{  } \;
\; {\rm rank} \frac{\partial f_i}{\partial x_{i+1}}
(t_1,x_1^{\ast},..., x_{i+1}^{\ast}) = m_i, \; \; \;  \;
i{=}1,...,\nu. \n{r4_4}$$ In addition, we introduce the
following notation
$$u^{\ast}{:=}x_{\nu{+}1}^{\ast}{\in}\rmnuplusl{=}\rmm; \;
\; \; x^{\ast}{:=}(x_1^{\ast},...,x_{\nu}^{\ast})^T {\in}
\rn {=} \rml{\times}...{\times}\rmnu   \n{r4_4a}  $$ To
prove theorem 1.1, it is sufficient to prove the following
statement.

{\bf Theorem 2.1.} {\sl Assume that, for system
(\ref{r4_1}), function $f$ has triangular form
      (\ref{r4_2}) and satisfies
     (I), (II). Choose
    $t_1 {\in} ]t_0,T[,$
 $z_{i}^{\ast}{\in} \rmi,$
    $x_i^{\ast} {\in} \rmi,$
 $i{=}1,...,\nu,$ and
$ x_{\nu{+}1}^{\ast}{=}u^{\ast}{\in}\rmm $  such that
(\ref{r4_4}), (\ref{r4_4a}) hold. Then there exist families
of controls
    $\{  {\overline u}_{x^T} (\cdot) {\}}_{x^T \in \rn}$ and
    $\{  {\tilde u}_{x^0} (\cdot) {\}}_{x^0 \in \rn}$ defined on
    $[t_1,T] $ and $[t_0,t_1]$ respectively
such that:

1) The maps   $x^T{\mapsto}{\overline u}_{x^T} (\cdot)$
    and  $x^0{\mapsto}{\tilde u}_{x^0} (\cdot)$
  are of classes $C(\rn; C([t_1,T]; \rmm))$
    and      $C(\rn; C([t_0,t_1]; \rmm))$ respectively.

     2) For each $x^0{\in}\rn,$ and each  $x^T{\in}\rn,$ we have:
      $  {\tilde u}_{x^0} (t_1){=}{\overline u}_{x^T} (t_1){=}u^{\ast}{=}x_{\nu{+}1}^{\ast}.$

     3) For each $x^0 {\in} \rn, $ and each  $x^T {\in} \rn,$  the trajectories
        $t \mapsto x(t,t_1, x^{\ast}, {\overline u}_{x^T} (\cdot))$    and
        $t \mapsto x(t,t_1, x^{\ast}, {\tilde u}_{x^0} (\cdot))$  are defined
     for all $t$ in   $[t_1,T]$ and  $[t_0, t_1]$ respectively, and
   $x(T,t_1,x^{\ast}, {\overline u}_{x^T} (\cdot)){=}x^T,$
       $x(t_0,t_1,x^{\ast}, {\tilde u}_{x^0} (\cdot)){=}x^0.$
    }

    Indeed, if theorem 2.1 is proven, then, the family
  $\{ u_{(x^0,x^T)} (\cdot) {\}}_{(x^0,x^T) \in \rn \times \rn}$
given by $$u_{(x^0,x^T)} (t) =  \left\{ \begin{array}{l}
     {\tilde u}_{x^0} (t) \; \; \; \; \; \; \; \;   \mbox{ if }\;  t \in [t_0,t_1[ \\
     {\overline u}_{x^T} (t) \; \; \; \; \; \; \; \;  \mbox{ if }\;  t \in [t_1,T],
   \end{array}\right. \; \; \; \; \;\quad\quad \mbox{ for all }  (x^0,x^T) {\in} \rn {\times} \rn  $$
satisfies theorem 1.1. On the other hand, to prove theorem
2.1, we need to construct only one of the two families of
controls (for instance $\{  {\overline u}_{x^T} (\cdot)
{\}}_{x^T \in \rn}$); the construction for the other  is
similar. The existence of the desired $\{  {\overline
u}_{x^T} (\cdot) {\}}_{x^T \in \rn}$ immediately follows
from the next statement (theorem 2.2), which is formulated
after the following notation.

Let  $p$ be in $\{ 1,...,\nu \}.$ By definition, put
 $$ k = m_1 + ...+m_p, \; \; y^{\ast} = (x_1^{\ast},...,x_p^{\ast}{)}^T \in \rk = \rml \times ...\times \rmp, $$
 $$ z^{\ast} = (z_1^{\ast},...,z_p^{\ast}{)}^T \in \rk, \; \; \; J:=[t_1,T], \n{r4_7}$$
and consider the following   $k$ -- dimensional control
system
      $$ \dot y  (t) = \varphi (t,y(t),v(t)), \; \; \; \; \; \; t \in J= [t_1,T], \n{r4_5}$$
where   $y{=} (x_1,...,x_p{)}^T {\in} \rk {=}
\rml{\times}...{\times}\rmp$ is the state,
      $v \in \rmpplusl$ is the control,
and
  $$\varphi(t,y,v) = \left( \begin{array}{l}
f_1(t,x_1,x_2) \\ f_2(t,x_1,x_2,x_3) \\
 \quad \ldots \ldots \ldots\quad \\
f_{p-1} (t,x_1,x_2,...,x_{p})\\ f_{p}
(t,x_1,x_2,...,x_{p},v)
\end{array}\right)\;\; \; \; \;   \n{r4_6}$$
for all  $(t,y,v){=}(t,x_1,...,x_p,v){\in}
  J {\times} \rmlplusdp {\times} \rmpplusl.$
  For each  $y{\in}\rk,$ each $\tau{\in}J,$   and each
   $v(\cdot) {\in} L_{\infty} (J; \rmpplusl),$ let
   $t{\mapsto} y (t,\tau,y, v(\cdot))$
 be the  trajectory, of system
    (\ref{r4_5}), that is defined by the control
   $v(\cdot)$ and by the initial condition
   $y(\tau,\tau,y, v(\cdot)) = y$ on some maximal subinterval
  $J_1 \subset J$ (where $\tau \in J_1$). If, for all $t \in J,$
we have $v(t)= v,$ where $v  \in \rmpplusl,$  then, we
denote this trajectory
 by $t \mapsto y(t,\tau, {y},v).$

 {\bf Theorem 2.2.}  {\sl    Assume that  $f$
has  triangular form (\ref{r4_2}) and satisfies
    (I), (II),   and let $t_1{\in}]t_0,T[,$
     $x^{\ast} {\in} \rn, $
     $x_{\nu+1}^{\ast}{=}u^{\ast} {\in} \rmm,$
and   $z_i^{\ast}{\in}{\rmi},$   $i{=}1,...,\nu,$  be such
that (\ref{r4_4}), (\ref{r4_4a}) hold. Choose any
$z_{\nu+1}^{\ast} {\in} \rmnuplusl.$ For some arbitrary
     $p {\in} \{ 1,...,\nu \},$ define
 $k,$    $y^{\ast} {\in} \rk,$  $z^{\ast} {\in} \rk,$ and $J$ by
     (\ref{r4_7}), and
consider system
     (\ref{r4_5}) with
     $\varphi$ given by
     (\ref{r4_6}).
Suppose that there exists a family
 $\{ y (\xi,\cdot) {=} (x_1(\xi,\cdot),...,x_p(\xi,\cdot){)}^T {\}}_{\xi \in \rk}$
  of functions
 of $J$ to $\rk$ such that:

     1)  For each $i{=}1,...,p,$  the map $\xi{\mapsto} x_i (\xi,\cdot)$   is of class
       $C(\rk ; C^1(J; \rmi)).$

     2)    For each   $\xi \in \rk,$
 we have:
  $ \dot x_i (\xi,t) = f_i(t,x_1(\xi,t),...,x_{i+1}(\xi,t))  $
whenever
 $1 {\leq} i {\leq} p{-}1,$  $i{\in}\nn,$   $t{\in}J.$

     3)  For each  $\xi \in \rk,$ we obtain:
     $y (\xi,t_1) {=} y^{\ast},$ $y(\xi,T) {=} \xi,$
$\dot x_p (\xi,t_1) {=} z_p^{\ast} {=}
f_p(t_1,y^\ast,x_{p{+}1}^{\ast}).$

Then, for system (\ref{r4_5}), there exists a family of
controls
     $\{ {\hat v}_{(\xi,\beta)} (\cdot)  {\}}_{(\xi,\beta) {\in} \rk {\times} \rmpplusl}$
such that:

      4)  The map  $(\xi,\beta) {\mapsto}  {\hat v}_{(\xi,\beta)} (\cdot)$ is of class
         $C(\rk {\times} \rmpplusl; C^1(J; \rmpplusl)).$

      5)   For each $ (\xi,\beta) {\in} \rk {\times} \rmpplusl,$ we obtain:
       ${\hat v}_{(\xi,\beta)} (T) {=} \beta,$
         ${\hat v}_{(\xi,\beta)} (t_1) {=} x_{p{+}1}^{\ast},$
         $\frac{d}{dt}{\hat v}_{(\xi,\beta)} (t_1) {=} z_{p{+}1}^{\ast}.$

      6) For each  $(\xi,\beta) \in \rk \times \rmpplusl,$ the trajectory
          $t {\mapsto} y (t,t_1,y^{\ast}, {\hat v}_{(\xi,\beta)} (\cdot))$
   is defined for all  $t{\in}J,$ and
          $ y (T,t_1,y^{\ast}, {\hat v}_{(\xi,\beta)} (\cdot)){=}\xi.$
       }

 Let us remark that, if a family of controls
  $\{ {\hat v}_{(\xi,\beta)} (\cdot)  {\}}_{(\xi,\beta) \in \rk \times \rmpplusl}$
satisfies conditions     4)   and  6) of theorem 2.2,
 then the map that assigns to each
     $(\xi,\beta){\in}\rk{\times}\rmpplusl$
the trajectory
     $t{\mapsto}y(t,t_1,y^{\ast}, {\hat v}_{(\xi,\beta)} (\cdot))$
is of class  $C(\rk {\times}\rmpplusl; C^1(J; \rk)).$
Therefore, theorem 2.2 means that, if, for $k$ --
dimensional system (\ref{r4_5}), there exists a family $\{
y(\xi,\cdot) {\}}_{\xi \in \rk}$ such that
      1)-3) hold, then the same is true for the
    $k+m_{p+1}$ -- dimensional dynamical extension, of (\ref{r4_5}),
that is determined by (\ref{r4_2}). For $ p{=} 1$ and
     $k {=} m_1,$ the construction of
     $ \{  y(\xi,\cdot)  { \}}_{\xi \in \rk}$ that satisfies
  1)-3)    is trivial:
 condition 2) consists of $p{-}1{=}0$ equalities; therefore,
$\{  y(\xi,\cdot)  { \}}_{\xi {\in} \rml}$
 given, for instance, by
     $$y(\xi,t) {=} x_1(\xi,t) {=} x_1^{\ast} {+} (t{-}t_1) z_1^{\ast} {+} \frac{(t{-}t_1{)}^2}{ (T{-}t_1{)}^2} (\xi {-} x_1^{\ast} {-} (T{-}t_1)z_1^{\ast}),
\; \; \;  \xi {\in} \rk {=} \rml, \; \; \; t {\in} J,$$
satisfies 1)-3). Using the induction on
 $p{=}1,...,\nu$ and theorem 2.2,
for $p {=} \nu,$ and $k{=}n,$ we get the existence of a
family of controls $\{  {\hat v}_{(\xi,\beta)} (\cdot)
{\}}_{(\xi,\beta){\in}\rn{\times}\rmm}$ such that
conditions 4)-6) of theorem 2.2 hold for $k{=}n;$ $m {=}
m_{\nu{+}1} {=} m_{p+1}.$ Fix an arbitrary $\beta {\in}
\rmm.$ Then, the family of controls
 $\{   {\overline u}_{x^T} (\cdot) {\}}_{x^T \in \rn}$
given by ${\overline u}_{x^T} (t){:=}{\hat v}_{(\xi,\beta)}
(t),$ $x^T{\in}\rn,$  $t{\in}J,$ satisfies the conditions
of theorem 2.1. The existence of the other family
 $\{   {\tilde u}_{x^0}    (\cdot) {\}}_{x^0 \in \rn}$
can be proved similarely. Thus, the main results of the
paper (theorems 1.1, 2.1, as well as 1.2 and 1.3) follow
directly from theorem 2.2 whose proof is the main objective
of our next efforts.

Throughout the paper, for each $r{>}0,$ and each  $y {\in}
\rk,$ by  $B_r (y)$ we denote the open ball $B_r (y){:=} \{
\eta {\in} \rk | \; | \eta {-} y| {<}r \}, $ and, for $A
\subset \rk,$ by $\overline{A}$ we denote the closure of
$A.$   In addition, we put $\zplus := \nn \bigcup \{ 0 \},$
and for  $r \in \rr, $ by  $[r]$ we denote $k \in \zz$ such
that $k {\leq} r < k{+}1.$

\vskip10mm
\begin{center}
 {\bf 3. Proof of theorem  2.2.   }
\end{center}

Let $p$ be in  $\{ 1,...,\nu \}.$ Assume that
 $\{   y(\xi,\cdot)  {\}}_{\xi \in \rk}$ satisfies conditions
 1)-3)   of theorem 2.2.

Let $x_i{=:}(x_{i,1},...,x_{i,m_i})^T {\in} \rmi$ and $f_i
(t,x,u) {=:} (f_{i}^1(t,x,u),...,f_i^{m_i}(t,x,u){)}^T$ be
the co\-ordi\-nate representations of vectors $x_i {\in}
\rmi $ and $f_i {\in} \rmi$ for  $i{=}1,...,\nu{+}1$ and
$i{=}1,...,\nu$ respectively. For every  $i{=}1,...,\nu,$
every $1{\leq}j_1<...<j_{m_i}{\leq}m_{i+1},$ and every
$(t,y,v) $ in $ [t_0,T]{\times}\rk{\times}\rmiplusl,$ by
$M_{j_1...j_{m_i}}^i (t,y,v)$ denote the matrix $\left(
\frac{\partial f_{i}^{\alpha}}{\partial x_{i{+}1,
j_{\beta}}} (t,y,v) \right)_{\alpha,\beta{=}1}^{m_i}$
(generated by the columns, of the matrix $\frac{\partial
f_i}{\partial x_{i{+}1}},$ numbered by
$j_1,j_2,...,j_{m_i}$). Using (\ref{r4_4}) and the implicit
function theorem, we get the existence of neighborhoods
  $E_{t_1, y^{\ast},x_{p+1}^{\ast}} {\subset} [t_0,T] {\times} \rk{\times}\rmpplusl$
and  $G_{t_1, y^{\ast},z_{p}^{\ast}} {\subset} [t_0,T]
{\times} \rk {\times}\rmp$
 of points   $(t_1, y^{\ast},x_{p{+}1}^{\ast})$
 and  $(t_1, y^{\ast},z_{p}^{\ast})$
  respectively and the existence of $p$ sequences of
indices  $1\leq j_1(i)<\ldots <j_{m_i} (i) \leq m_{i+1},$
    $i \in \{ 1,...,p \},$ such that, first,
 $$ \forall (t,y,v) \in E_{t_1,y^{\ast},x_{p+1}^{\ast}}\; \; \;  \forall i
\in \{ 1,...,p\}\; \; \; \; {\rm det}
M_{j_1(i),...,j_{m_i}(i)}^i (t,y,v) \not= 0, \n{r4_9}$$
and, second, there exists a map
 $(t,y,z_{p+1}) \mapsto \phi(t,y,z_{p+1})$  of class
  $C(G_{t_1, y^{\ast},z_{p}^{\ast}};\rmpplusl)$
such that $\frac{\partial \phi}{\partial y}$ and
$\frac{\partial \phi}{\partial z_{p{+}1}}$ are continuous
on $ G_{t_1, y^{\ast},z_{p}^{\ast}},$ and
  $$ \phi (G_{t_1, y^{\ast},z_{p}^{\ast}}) {\subset}
\left\{ (x_{p{+}1, 1}, \ldots, x_{p{+}1, m_{p{+}1}}{)}^T
{\in} \rmpplusl |\;\;\; \;  \forall j {\in} \{ 1,...,
m_{p{+}1} \} \right. $$ $$ \left. \left( j \notin \{
j_1(p),...,j_{m_p}(p)\} \right)  {\Rightarrow} \left(
x_{p{+}1,j } = x_{p{+}1, j}^{\ast}\right) \right\}; \; \;
\; \phi(t_1,y^{\ast},z_p^{\ast}){=}x_{p+1}^{\ast};$$
      $$\forall (t,y,z_p) \in G_{t_1, y^{\ast},z_{p}^{\ast}}\; \; \;
\; \;  \; \; \; \; \;  \left(t,y,\phi(t,y,z_p)\right) \in
E_{t_1, y^{\ast},x_{p+1}^{\ast}} ; \n{r4_10}$$
      $$\forall (t,y,z_p) \in G_{t_1, y^{\ast},z_{p}^{\ast}}\;
 \; \; \; \;  \quad\quad
\; \; \; \; \; f_p (t,y,\phi(t,y,z_p)) =z_p.  \n{r4_11}$$
Let us prove the existence of $\sigma (\cdot) {\in} C(\rk;
]0, T{-}t_1[)$ such that
     $$\forall \xi \in \rk\;  \; \;  \forall s \in [t_1, t_1+ \sigma(\xi)] \;
\;\;\;\;      \left(s,y(\xi,s),\dot{x_p} (\xi,s) \right)
\in G_{t_1, y^{\ast},z_{p}^{\ast}}. \n{r4_12}$$ For this,
it suffices to prove the existence of $\{ {\delta}_l
{\}}_{l=1}^{\infty} \subset ]0,T{-}t_1[$ such that, for
each
        $l {\in} \nn,$  each
       $\xi {\in} \rk$ such that
       $|\xi| {\leq} l,$ and each
      $s {\in} [t_1, t_1{+} {\delta}_l],$ we have:
     $(s,y(\xi,s),\dot{x_p}  (\xi,s)) \in   G_{t_1, y^{\ast},z_{p}^{\ast}}.$
Then, without loss of generality, we may assume that
 $0 < {\delta}_{l{+}1} {\leq} {\delta}_{l} {<} T{-}t_1$
for all  $l {\in} \nn, $ and then the continuous function
$\sigma(\cdot)$ given by
   $$ \sigma(\xi) = \left\{ \begin{array}{l}
    {\delta}_{[ |\xi| +1]}  + ({\delta}_{[|\xi| +1] +1} - {\delta}_{[|\xi| +1]}) (|\xi| - [|\xi|]) \; \; \; \; \;    \mbox{ if } |\xi| \notin \zplus  \\
   {\delta}_{|\xi| +1} \; \; \; \; \; \; \quad\quad\quad\quad \quad\quad\quad\quad\quad\quad\quad\quad\quad\quad\mbox{ if } |\xi| \in \zplus,
 \end{array}\right.$$
  satisfies (\ref{r4_12}).
Assume the converse (i.e. that such $\{ \delta_l \}$ does
not exist). Then, there exist $l_0 {\in} \nn,$ a sequence
 $\{ {\overline \xi}_q  {\}}_{q{=}1}^{\infty} {\subset} \overline{B_{l_0} (0)}{\subset}\rk,$
and a sequence
  $\{ {\overline \delta}_q  {\}}_{q=1}^{\infty} \subset ]0,T{-}t_1[$
such that
  ${\overline \delta}_q {\rightarrow} 0$ as $q {\rightarrow}  + \infty,$
and for each
  $ q {\in} \nn$  we have:
  $$\{ (s,y({\overline \xi}_q,s), \dot{ x_p}({\overline \xi}_q,s)) {\in}
 J {\times} \rk {\times} \rmp |\; \;  s {\in} [t_1, t_1 {+} {\overline \delta}_q]
\} \not\subset G_{t_1, y^{\ast},z_{p}^{\ast}}.$$ Since
$|{\overline \xi}_q | \leq l_0 $ for all  $q \in \nn,$ we
obtain that there exist
  ${\overline \xi} {\in} \rk$
  ($|{\overline \xi}| {\leq} l_0$) and a subsequence
   $\{ {\overline \xi}_{q_{\Upsilon}} {\}}_{\Upsilon {=}1}^{\infty}$
of  sequence
   $\{   {\overline \xi}_q {\}}_{q=1}^{\infty}$
such that ${\overline \xi}_{q_{\Upsilon}} {\to} {\overline
\xi}$ as $\Upsilon {\to} \infty.$ Without loss of
generality, we may assume that  ${\overline \xi}_q
{\rightarrow} {\overline \xi} $ as
   $q {\rightarrow} 0.$
 For
    ${\overline \xi},$
we obtain from conditions 1) and 3) of theorem 2.2 that
there exists
    $\delta ({\overline \xi}) {\in} ]0,T{-}t_1[$
such that, for each
     $s {\in} [t_1, t_1 {+} \delta ({\overline \xi})],$
we have
    $(s,y({\overline \xi},s), \dot{ x_p}({\overline \xi},s)) {\in}
G_{t_1, y^{\ast},z_{p}^{\ast}}.$ Then, from condition 1) of
theorem 2.2, we get the existence of  $\rho {>} 0$ such
that, for each
       $\xi {\in} B_{\rho} ({\overline \xi}),$         the inclusion
$(s,y({ \xi},s), \dot{ x_p}({ \xi},s)) {\in}   G_{t_1,
y^{\ast},z_{p}^{\ast}}$
 holds for all
$s {\in} [t_1, t_1 {+} \delta({\overline \xi})]$ as well.
In particular, this is true for   $\xi = {\overline \xi}_q$
 whenever   $q\geq q_0 $
(where    $q_0 {\in} \nn$ is any number
 such that for all  $q{\geq} q_0$ we get
    ${\overline{\xi}}_q  {\in} B_{\rho} ({\overline \xi})$).
This contradicts the choice of
   $\{   {\overline \delta}_q {\}}_{q=1}^{\infty}$
and      $\{   {\overline \xi}_q {\}}_{q=1}^{\infty}$ and
proves the existence of    $\sigma(\cdot) {\in} C(\rk;
]0,T{-}t_1[)$
 such that (\ref{r4_12}) holds.

Consider the family of $k$ - dimensional control systems
     $$\dot z (t){=} \frac{\partial \varphi}{\partial y} ( t, y(\xi,t),
\phi(t,y(\xi,t), \dot{x_p} (\xi,t))) z(t) {+}
 \frac{\partial \varphi}{\partial v} ( t, y(\xi,t), \phi(t,y(\xi,t),
\dot{x_p} (\xi,t))) w(t) , \;    $$ $$  \; \;  t {\in}
[t_1, t_1{+}\sigma(\xi)] \n{r4_13}$$ where      $z {\in}
\rk$ is the state,   $w {\in} \rmpplusl$ is the control,
and    $\xi \in \rk $ is the parameter of the family. From
(\ref{r4_9}), (\ref{r4_10}), (\ref{r4_12}), and from
 lemma 4.1, we get the existence of
$k$ families of controls $\{  w_j (\xi, \cdot) {\}}_{\xi
{\in} \rk},$ $ j {=} 1,...,k,$ such that
 $$   w_j (\cdot,\cdot) \in C(\{(\xi,t) {\in} \rk {\times}J |\; \; t_1{\leq} t {\leq} t_1 {+} \sigma(\xi)\} ; \rmpplusl),$$
 $$  \frac{\partial w_j}{\partial t} (\cdot,\cdot) \in C(\{(\xi,t) {\in} \rk {\times} J |\; \; t_1{\leq }t {\leq} t_1 {+} \sigma(\xi)\} ; \rmpplusl), \; \; \;
 j {=} 1,...,k, \n{r4_19}$$
and such that for each $ j {=} 1,...,k,$ and each
  $\xi \in \rk,$
the control $ w_j(\xi,\cdot)$
 is defined on
  $[t_1,t_1 {+} \sigma(\xi)],$
steers
 $0 {\in} \rk$ into
   $e_j {=} (0,...,0,1,0,...,0{)}^T {\in} \rk$
   (the $j$ -- th unit vector of $\rk$)
in time     $[t_1, t_1 {+} \sigma(\xi)]$ with respect to
    (\ref{r4_13}), and satisfies the boundary conditions
$$ {w}_j(\xi,t_1) {=} 0 {\in} \rmpplusl, \; \;   {\dot w}_j
(\xi,t_1) {=} 0 {\in} \rmpplusl,\;\;   {w}_j (\xi,t_1 {+}
\sigma(\xi)) {=} 0 {\in} \rmpplusl, $$ $$  {\dot w}_j
(\xi,t_1 {+} \sigma(\xi)) {=} 0 {\in} \rmpplusl.
\n{r4_18}$$ For each  $\lambda {=}
({\lambda}_1,...,{\lambda}_k{)}^T {\in} \rk,$ and each
$\xi {\in} \rk,$ let
   $v_{\lambda} (\xi,\cdot)$ be the control
defined on
   $[t_1, t_1 {+} \sigma (\xi)]$ by
$$v_{\lambda}(\xi, t) {=} \phi(t,y(\xi,t), {\dot x}_p
(\xi,t)) {+} \sum\limits_{j=1}^{k} {\lambda}_j w_j (\xi,t),
\; \; \; \mbox{ for all }\;  t {\in} [t_1, t_1 {+} \sigma
(\xi)]. \n{r4_20}$$ For each $\xi {\in} \rk,$ and each
$\lambda {\in} \rk$ such that
    $t {\mapsto} y(t,t_1, y^{\ast} , v_{\lambda} (\xi,\cdot))$
is defined for all
    $t{\in}[t_1,t_1{+} \sigma(\xi)],$ we put, by definition,
    $y_{\lambda}(\xi,t){:=}y(t,t_1,y^{\ast}, v_{\lambda} (\xi,\cdot))$
whenever  $t {\in} [t_1, t_1{+}\sigma(\xi)], $ and, then,
for each   $\mu {=} ({\mu}_1,...,{\mu}_k{)}^T {\in} \rk,$
by $z_{\mu, \lambda} (\xi,\cdot)$ we denote the trajectory,
of the system
  $$   {\dot z} (t) {=} \frac{\partial \varphi}{\partial y} (t, y_{\lambda}
(\xi,t) , v_{\lambda} (\xi,t)) z(t) {+} \frac{\partial
\varphi}{\partial v} (t, y_{\lambda} (\xi,t) , v_{\lambda}
(\xi,t)) w(t), \; \; \; \;
   t {\in} [t_1, t_1 {+} \sigma(\xi)], \n{r4_21}$$
 with states  $z {\in} \rk,$ and controls
       $w {\in} \rmpplusl,$ that is defined on
    $[t_1, t_1 {+} \sigma (\xi)]$ by the control
    $w_{\mu} (\xi,\cdot) {=} \sum\limits_{j=1}^{k} {\mu_j}w_j(\xi,\cdot),$
and by the initial condition
     $z_{\mu,\lambda}(\xi,t_1){=}0 {\in} \rk.$

Define the families   $\{ \Phi(\xi,\cdot) {\}}_{\xi {\in}
\rk}$ and     $\{ \Psi(\xi,\cdot,\cdot) {\}}_{\xi {\in}
\rk}$ of maps from  $\rk $ and
  $\rk {\times}\rk$ respectively
to  $\rk$ as follows: for each
     $\xi {\in} \rk,$ each
   $ \mu {\in} \rk,$  and each
     $ \lambda {\in} \rk$ such that
the trajectory     $ t {\mapsto} y(t,t_1,y^{\ast},
v_{\lambda}(\xi,\cdot))$ is defined for all
    $t {\in} [t_1, t_1 {+} \sigma (\xi)],$
by definition, put:
     $\Phi (\xi,\lambda) {=} y_{\lambda}(\xi, t_1{+}\sigma(\xi)),$
  $\Psi(\xi,\mu,\lambda) {=} z_{\mu,\lambda}(\xi,t_1{+}\sigma(\xi)).$
(Note that, by the construction, from (\ref{r4_11}),
(\ref{r4_12}), (\ref{r4_20}), and from condition 2) of
theorem 2.2, we get
  $ y_{\lambda}(\xi, t){|}_{\lambda{=}0} {=} y(\xi,t),$
$t{\in} [t_1,t_1{+}\sigma(\xi)],$ and, therefore,
  $\Phi (\xi,0) {=} y(\xi, t_1{+}\sigma(\xi))$
for all $\xi{\in}\rk$).

 {\bf  Lemma 3.1.}  {\sl   1)  There exists
a function $\varepsilon (\cdot) {\in} C(\rk; ]0,+\infty[)$
            such that, for each  $(\xi,\lambda)  $ in
$\Pi {:=} \{ (\xi,\lambda) {\in} \rk {\times}\rk |\;
\lambda {\in} B_{\varepsilon (\xi)} (0)\},$ the trajectory
$t {\mapsto} y(t,t_1 , y^{\ast}, v_{\lambda} (\xi,\cdot))$
is defined for all
           $t{\in}[t_1,t_1 {+} \sigma(\xi)],$
and, therefore,  $\Phi(\xi,\lambda)$ and $\Psi
(\xi,\cdot,\lambda)$ are well defined.

  2) For each $\xi {\in} \rk, $
the map  $\lambda {\mapsto} \Phi (\xi, \lambda)$ is
differentiable for all  $ \lambda {\in}
B_{\varepsilon(\xi)} (0),$ and
 $$ \frac{\partial \Phi}{\partial \lambda} (\xi,\lambda) \mu = \Psi (\xi,\mu, \lambda) \; \;
\; \; \mbox{ whenever } \; \; \lambda{\in}B_{\varepsilon
(\xi)} (0), \; \;
          \mu{\in}\rk. $$

    3)  The maps    $(\xi,\lambda) {\mapsto} \Phi(\xi,\lambda)$ and
     $(\xi,\lambda) {\mapsto} \frac{\partial \Phi}{\partial \lambda} (\xi,\lambda)$
are of classes     $C(\Pi; \rk)$ and $C(\Pi; \rktk)$
respectively.
        }

Lemma 3.1 immediately follows from the well-known theorems
on the differentiability of the solution of the Cauchy
problem w.r.t. a parameter.

Using the definition of $\{w_j(\xi,\cdot){\}}_{\xi \in
\rk}$ and lemma 3.1, for each $\xi {\in} \rk,$ and each
$j{=}1,...,k,$ we get: $\Psi (\xi,e_j,\lambda) {|}_{\lambda
{=} 0} {=} e_j,$
 i.e., $\frac{\partial \Phi}{\partial \lambda} (\xi,0){=}I,$
where $I{\in}\rktk$ is the identity matrix. Fix some $\rho
{>} 0$ such that each matrix $A {\in} \rktk$ that satisfies
the inequality $\parallel A{-}I{\parallel} {<} 2 \rho$ is
positive definite.

{\bf Lemma 3.2.}  {\sl  There exist  functions
 $ \varepsilon_1(\cdot)$ and  $\varepsilon_2 (\cdot),$
of class
   $C(\rk; ]0,+\infty[),$ such that for each
  $\xi {\in} \rk$ we have:
$$   \varepsilon_1(\xi) < \frac{1}{2}  \varepsilon (\xi)
\quad\quad\quad\quad\quad\quad \n{r4_22}$$ $$\forall
\lambda \in \overline{B_{\varepsilon_1 (\xi)} (0) }\quad\;
\;
\parallel  \frac{\partial \Phi}{\partial \lambda} (\xi,\lambda) - I
\parallel < \rho \n{r4_23}$$
$$ \overline{B_{\varepsilon_2  (\xi) }
(y(\xi,t_1+\sigma(\xi)))} \subset \Phi (\xi,
B_{\varepsilon_1  (\xi)} (0)). \n{r4_24}$$ }

The proof of lemma 3.2, which is omitted, is the same as
the proof of lemma 3.4 from \cite{ssp_op1}. The
construction of $\varepsilon_1 (\cdot)$ and $\varepsilon_2
(\cdot)$ is similar to that of the continuous function
$\sigma (\cdot)$ introduced above.

Along with system (\ref{r4_5}), we consider the following
$k$ - dimensional control system $$  \left\{
\begin{array}{l} \dot x_i (t)  = f_i (t,x_1 (t),
...,x_{i+1}(t)), \; \; \; \; \; 1 {\leq} i {\leq} p{-}1,
\;\; i {\in} \nn; \\ \dot x_p (t) = \omega (t),
\end{array}\right. \; \; \; \; \; \; \; \;     t \in J \n{r4_25}  $$
with states $y {=} (x_1,...,x_p{)}^T {\in} \rmlplusmp {=}
\rk$ and controls $\omega {\in} \rmp.$

 For each $y {\in} \rk,$
each $\tau {\in} J,$ and each $\omega(\cdot) {\in}
L_{\infty}(J; \rmp),$ by $t {\mapsto} z(t,\tau,y, \omega
(\cdot))$ we denote the trajectory, of (\ref{r4_25}), that
is defined by the control $\omega(\cdot)$ and by the
initial condition $z(\tau,\tau,y,\omega(\cdot)){=}y$  on
some maximal subinterval
 $J_1 {\subset} J$   $(\tau {\in} J_1).$

Let us remark that, from conditions 1)-3) of theorem  2.2,
and from the definition of $\{ \Phi (\xi,\cdot) {\}}_{\xi
\in \rk},$ we obtain: $$ \forall {\xi} {\in} \rk  \;\;
\forall t {\in}J \;\;   y(\xi,t) {=} z(t,T,\xi, \dot x_p
(\xi,\cdot));  \; \; \;\;  \mbox{ and } \; \; \; \Phi
(\xi,0) {=} z(t_1{+}\sigma(\xi),T,\xi, \dot x_p
(\xi,\cdot)). \n{r4_26}$$

{\bf Lemma 3.3.} {\sl  There exists
 $\delta(\cdot) \in C(\rk ; ]0, + \infty[)$ such that, for each
 $\xi \in \rk,$ and each
$\omega(\cdot) \in L_{\infty} (J; \rmp)$ that satisfies $
|| \omega(\cdot) - \dot x_p  (\xi,\cdot) {||}_{L_{\infty}
(J; \rmp) } <  \delta (\xi),$ the trajectory $t {\mapsto}
z(t,T,\xi,\omega (\cdot))$ is defined for all $t{\in}J,$
and satisfies the conditions $$\forall t \in J \;\quad \;
|z(t,T,\xi,\omega(\cdot)) {-} y(\xi,t)| <
\frac{\varepsilon_2 (\xi)}{ 4} ; \n{r4_27}$$ $$| z(t_1 +
\sigma(\xi),T,\xi,\omega(\cdot)) {-} \Phi (\xi, 0)| <
\frac{\varepsilon_2 (\xi)}{4}. \n{r4_28}$$ }

{\bf Proof of lemma 3.3.}   Like in the construction of
$\sigma (\cdot),$ it suffices to prove the existence of $\{
{\overline{\delta}}_r   {\}}_{r=1}^{\infty}$
($0<{\overline{\delta}}_{r+1} < {\overline{\delta}}_r$ for
all $r {\in} \nn$) such that, for each
 $r {\in} \nn, $ each
$\xi {\in} \overline{B_{r}(0)},$ and each $\omega(\cdot)
{\in} L_{\infty} (J;\rmp),$ the inequality ${||} \omega
(\cdot) {-} \dot x_p (\xi,\cdot) {||}_{L_{\infty} (J
;\rmp)}    < {\overline{\delta}}_r$
 implies that $t {\mapsto} z(t,T,\xi,\omega(\cdot))$
 is defined for all $t{\in}J, $  and
$$ \forall t {\in} J \; \;  |z(t,T,\xi,\omega(\cdot)) {-}
z(t,T,\xi,\dot x_p   (\xi,\cdot)) | < \frac{1}{4}
\min\limits_{\xi \in  \overline{B_r (0)}} \varepsilon_2
(\xi).$$ The existence of such $\{   {\overline{\delta}}_r
{\}}_{r=1}^{\infty}$ follows from standard arguments based
on the Gronwall-Bellman lemma. The proof of lemma 3.3 is
complete.

{\bf   Lemma 3.4.}  {\sl  Let family
  $\{ y(\xi,\cdot)  {\}}_{\xi \in \rk}$ be such that
conditions  1)-3)    of theorem 2.2 hold. Then, for system
(\ref{r4_5}), there exist a function   $M(\cdot) \in C(\rk;
]0,+ \infty[)$ and a family   $\{  v(\xi,\cdot) {\}}_{\xi
\in \rk}$
 of controls
  defined on  $J$  such that the following conditions hold:

1)  For each $\xi {\in} \rk,$ the control
    $v(\xi, \cdot)$  is a piecewise constant function
 on   $J,$  and the map      $\xi {\mapsto} v(\xi, \cdot) $
 is of class    $C(\rk; L_1(J ; \rmpplusl)).$

2) For each $ \xi {\in} \rk, $ the trajectory
    $t {\mapsto} y(t,T,\xi,v(\xi,\cdot))$ is defined for all
      $t{\in}J,$ and
    $$\forall t \in J \quad\; \; \;  \; \; |\dot x_p (\xi,t) - f_p(t,y(t,T,\xi,v(\xi,\cdot)),v(\xi,t))| < \delta(\xi).$$

3)  For each $\xi {\in} \rk,$    we have:
  $ \parallel v(\xi,\cdot) {\parallel}_{L_{\infty} (J ; \rmpplusl)}  {\leq} M(\xi).$

}

\vskip10mm
\begin{center}
{\bf 3.1.  Proof of lemma 3.4.}
\end{center}

Let $\{   R_q   {\}}_{q=1}^{\infty} \subset \nn$ be an
arbitrary sequence such that
 $R_1 {=} 1,$
$R_{q{+}1} {>} R_q {+}1,$ $q {\in} \nn.$  Let us recall
that for each $ \eta {\in} \rk,$  and each $R{>}0$ by $B_R
(\eta)$ we denote the set $\{  \zeta {\in} \rk |\; | \zeta
{-} \eta| {<} R \}.$ By definition, put $$ \delta_q {=}
\frac{1}{2} \min\limits_{\xi \in \overline{ B_{R_{q{+}1}}
(0)} } \delta (\xi), $$ $$ M_q {=} \max\limits_{\xi \in
\overline{B_{R_q}(0) }} \parallel y(\xi,\cdot)
{\parallel}_{C(J;\rk)}  + \max\limits_{\xi \in
\overline{B_{R_q}(0) }} \varepsilon_2 (\xi) {+}1, \; \;\; q
{\in} \nn ; \n{r42_29}$$ $$  K_q = \{ y \in \rk | \; |y|
\leq M_q \} \subset \rk;  \; \; \; d_q = M_{q+2}+1, \; \;
\; q \in \nn ; \n{r42_30}$$ $$ W_q = \{ z_p \in \rmp | \;
|z_p| \leq \max\limits_{\xi \in \overline{ B_{R_{q}} (0)}
}\parallel \dot x_p (\xi,\cdot)  {\parallel}_{C(J; \rmp)}
+1 \} \subset \rmp , \; \;
 q {\in} \nn ; \; \;\; \; \; \; \; \; \; \n{r42_31}$$
$$\Xi_1 = \overline{B_{R_1} (0)}; \; \; \Xi_{q+1} =
\overline{ B_{R_{q+1}}(0)}\setminus \overline{
B_{R_{q}}(0)}, \; \; q \in \nn; \n{r42_32}$$ $$ E_1 =
\overline{B_{R_1} (0)} {\times} J {\times} K_1; \; \; \;
E_{q+1} = E_q {\bigcup} \left( \left(
\overline{B_{R_{q+1}}(0)} \setminus B_{R_q} (0) \right)
 \times J \times K_{q+1} \right), \; \;   q {\in} \nn ; \n{r42_33} $$
$$ E = \bigcup\limits_{q=1}^{\infty} E_q . \n{r42_34} $$

Fix an arbitrary $q {\in} \nn.$ For each $N {\in} \nn,$ by
definition, put: $${\Lambda}_N^q = \left\{  (t,y,z_p) {\in}
J {\times} K_{q{+}1} {\times} \rmp | \;\;  \exists \;
{\overline v} {\in} \rmpplusl  \; \;\left( |{\overline v}|
{\leq} N \right) \wedge
 ( |z_p {-} f_p(t,y,{\overline v})|  {<} \frac{\delta_q}{3} ) \right\}.$$
Since $f_p$ is continuous,
 each
${\Lambda}_N^q$  $(N \in \nn)$ is an open set of the metric
space $J {\times} K_{q{+}1} {\times} \rmp$ equipped with
the usual metric generated by the norm of $\rr {\times} \rk
{\times} \rmp.$
 In addition, $J {\times} K_{q{+}1} {\times} W_q$
is a compact set of this metric space, and ${\Lambda}_N^q
{\subset} {\Lambda}_{N+1}^q$ for all $N {\in} \nn.$  Then,
since $J \times K_{q+1} \times W_q \subset
\bigcup\limits_{N=1}^{\infty} {\Lambda}_{N}^q,$ from (II),
we get the existence of   $N_{0} (q) {\in} \nn$ such that
$J {\times} K_{q{+}1} {\times} W_q {\subset} {\Lambda}_{N_0
(q)}^q.$ By definition, put: $$U_q = \{ v {\in} \rmpplusl |
\; \; |v| \leq N_0 (q) \}. \n{r42_35}$$ Thus, for each $q
{\in} \nn,$ we get the existence of
 the compact set $U_q$ defined by (\ref{r42_35})
such that for every $(t,y,z_p) {\in} J {\times} K_{q{+}1}
{\times} W_q$ there exists $v {\in} U_q $ such that $|z_p -
f_p (t,y,v)| < \frac{{\delta}_q}{3}.$

In addition, for each $q {\in} \nn,$ by definition, put
$$L_q = \frac{1}{2}\; \;  {( \max\limits_{\begin{array}{c}
t \in J \\ y \in \overline{B_{d_q} (0)} \\ v \in U_{q+2}
\end{array}}  |\varphi (t,y,v)| +1)}^{-1}, \; \; \; q \in \nn, \n{r42_36} $$
  and let $L(\cdot) \in C(\rk ; ]0, + \infty[)$
be an arbitrary function such that $$L_{q+1} \leq L(\xi)
\leq L_q  \; \; \; \; \;  \mbox{  whenever } \xi \in
{\Xi}_q, \; \; q \in \nn \n{r42_36_a}$$ (for instance,
$L(\cdot)$ given by $$ L(\xi)  = \left\{ \begin{array}{l}
L_2 \; \; \; \; \;  \; \; \; \; \; \; \; \; \; \; \; \; \;
\; \; \; \; \; \; \; \; \; \; \; \; \; \; \; \; \; \;
\mbox{ if } \xi \in {\Xi}_1  \\ \frac{R_q - |\xi|}{R_q
-R_{q-1}}L_q + \frac{|\xi| - R_{q-1}}{R_q - R_{q-1}}
L_{q+1}  \; \;  \mbox{ if }  \xi \in {\Xi}_q,\;   q \geq 2,
\;  q \in \nn
\end{array}\right.  $$
is continuous and satisfies (\ref{r42_36_a})).

For each $(\xi,t,y) {\in} E,$ let $q {\in} \nn$ be such
that $\xi {\in} {\Xi}_q.$ Then $y {\in} K_{q+1},$ and, by
the construction of $U_q,$ there exists $v_{\xi,t,y} {\in}
U_q$ such that $|\dot x_p (\xi,t) - f_p (t,y,v_{\xi,t,y})|
< \frac{{\delta}_q}{3}.$ Since $f_p$  is continuous, from
condition   2)   of theorem 2.2, we get the existence of an
interval $I_{\xi,t,y} = ]t {-} {\overline
\theta}_{\xi,t,y}, \; t{+} {\overline \theta}_{\xi,t,y}[,$
where ${\overline \theta}_{\xi,t,y} > 0,$ and ${\overline
\sigma} (\xi,t,y), $ ${\overline \rho} (\xi,t,y) $ are
numbers from  $ ]0, \frac{1}{2}[,$ such that, for each
$\eta {\in} B_{{\overline \sigma} (\xi,t,y)} (\xi),$ each
$s {\in} I_{\xi,t,y} {\bigcap} J,$  and each $ z {\in}
B_{{\overline \rho} (\xi,t,y)} (y),$ we have: $$|\dot x_p
(\eta,s) - f_p (s,z,v_{\xi,t,y})| < \delta_q.  \n{r42_37}$$
 Choose an arbitrary
 ${\overline \tau}_{\xi,t,y} \in ]0, {\overline \theta}_{\xi,t,y} [$
such that $${\overline \tau}_{\xi,t,y} < \min \left\{
L_{q+2} \; {\overline \rho} (\xi,t,y)   ,\; \; T-t_1
\right\}   \; \; \; \; \; \; \; \; \; \; \; \; (\xi \in
{\Xi}_q). \n{r42_38}$$ Let ${\theta}_{\xi,t,y}
(\cdot,\cdot)$ and
 ${\tau}_{\xi,t,y} (\cdot,\cdot),$
be the functions of $\rk {\times} \rk $
 to $\rr$ and
$T_{\xi,t,y} \subset \rk {\times}\rr {\times}\rk$ be the
open set that are given by $${\theta}_{\xi,t,y} (\eta,z) =
t{+} {\overline \tau}_{\xi,t,y} {-} L_{q{+}2}  |z{-}y|  {-}
4(T{-}t_1) \max \left\{ 0, \; \frac{|\eta
{-}\xi|}{{\overline \sigma}(\xi,t,y) } {-}
\frac{1}{2}\right\}, \; \;
 \xi {\in} {\Xi}_q;  \n{r42_39}$$
$${\tau}_{\xi,t,y} (\eta,z) {=}  t{-} {\overline
\tau}_{\xi,t,y} {+} L_{q{+}2}  |z{-}y|  {+} 4(T{-}t_1) \max
\left\{ 0, \; \frac{|\eta {-}\xi|}{{\overline
\sigma}(\xi,t,y) } {-} \frac{1}{2} \right\}, \; \;
 \xi {\in} {\Xi}_q;  \n{r42_40}$$
$$ T_{\xi,t,y} = \{ (\eta,s,z) {\in} \rk {\times} \rr
{\times} \rk | \; {\tau}_{\xi,t,y} (\eta,z) {<} s {<}
{\theta}_{\xi,t,y}(\eta,z) \}.   \n{r42_41}$$ From
(\ref{r42_38})  it follows that $T_{\xi,t,y} \subset
B_{{\overline \sigma} (\xi,t,y)} (\xi) {\times} I_{\xi,t,y}
{\times} B_{{\overline \rho} (\xi,t,y)} (y),$ and,
therefore, for each $(\eta,s,z) {\in} T_{\xi,t,y}$
 such that $s {\in} J, $ we obtain (\ref{r42_37}); combining this
with the inequalities ${\overline \sigma} (\xi,t,y) {<}
\frac{1}{2},$
 and ${\overline \rho} (\xi,t,y) {<} \frac{1}{2},$ we get
$$ \forall (\eta,t,z) \in T_{\xi,t,y} \; \; \; \; \; \left(
|\eta-\xi| < \frac{1}{2} \right) \wedge \left( |z-y| <
\frac{1}{2} \right). \n{r42_42}$$
 Then, for each $(\eta,s,z) {\in} T_{\xi,t,y},$ we get:
$\eta \in \overline{B_{R_{q+1}} (0)} .$  By the definition
of ${\delta}_q, $ this yeilds: $$ \forall (\eta,s,z) \in
T_{\xi,t,y} \bigcap (\rk {\times} J {\times}\rk) \; \; \;
\; \;\; \; \; \;  | \dot x_p (\eta,s) -  f_p
(s,z,v_{\xi,t,y}) | < \delta (\eta).   \n{r42_43}$$ Thus,
we have constructed  the family of pairs $\{ (T_{\xi,t,y},
v_{\xi,t,y}) {\}}_{(\xi,t,y) \in E},$ where each
 $T_{\xi,t,y}$ is an open subset of $\rk {\times}\rr{\times}\rk$ such that
  (\ref{r42_43}) holds.

For the open covering $\{ T_{\xi,t,y} {\}}_{(\xi,t,y) \in
E}$ of $E,$ choose a countable locally finite subcovering
$\{ T_{{\xi}_r,t_r,y_r} {\}}_{r=1}^{\infty}$ of $E$ as
follows.
 For $q=1,$ let
$\{ T_{{\xi}_r,t_r,y_r} {\}}_{r=1}^{r_q = r_1}$ be a finite
subcovering
 selected from the open covering $\{ T_{{\xi},t,y} {\}}_{(\xi,t,y) \in E_1 }$
 of the compact set $E_1$
(in particular, we have
 $(\xi_r,t_r,y_r) \in E_1,$ $r =1,...,r_1$).
Assume that, using the induction over $q {\in} \nn,$
 for some $q {\in} \nn,$ we have already constructed
 a finite covering $\{ T_{{\xi}_r,t_r,y_r} {\}}_{r=1}^{r_q},$
 of the compact set $E_q,$ such that
  $(\xi_r,t_r,y_r) {\in} E_q,$
$r {=} 1,...,r_q.$ Then, from (\ref{r42_42}),  and from
(\ref{r42_33}), we obtain that the compact set $E_{q+1}
\setminus (\bigcup\limits_{r=1}^{r_q} T_{\xi_r,t_r,y_r})$
is not empty. Select from its open covering $\{
T_{{\xi},t,y} {\}}_{(\xi,t,y) \in E_{q+1} \setminus
(\bigcup\limits_{r=1}^{r_q} T_{\xi_r,t_r,y_r}) }$ a finite
subcovering $\{ T_{{\xi}_r,t_r,y_r} {\}}_{r=r_q +
1}^{r_{q+1}},$ and obtain the finite covering $\{
T_{{\xi}_r,t_r,y_r} {\}}_{r= 1}^{r_{q+1}},$ of the compact
set $E_{q{+}1},$ such that $({\xi}_r,t_r,y_r) {\in}
E_{q{+}1},$ $r{=}1,...,r_{q{+}1}.$

By the construction, for $r{\geq} r_{q{+}1} {+}1$ $(q {\in}
\nn)$  we have $({\xi}_r,t_r,y_r) \in
\bigcup\limits_{m{=}1}^{\infty}\left(E_{q+m+1} \setminus
E_{q+m}\right);$ then, from  (\ref{r42_42}) we obtain
$T_{{\xi}_r,t_r,y_r} \bigcap \left( \overline{B_{R_q+
\frac{1}{2}} (0)} \times \rr\times\rk \right)  =
\emptyset,$ and, therefore, $$T_{{\xi}_r,t_r,y_r} \bigcap
\left(\bigcup\limits_{j=1}^{r_q} T_{{\xi}_j,t_j,y_j}\right)
= \emptyset \; \;\; \; \;   \mbox{ whenever } r \geq
r_{q+1}+1, \; \; q \in \nn. \n{r42_44}$$ In particular,
from this, we obtain that the covering $\{
T_{{\xi}_r,t_r,y_r} {\}}_{r= 1}^{\infty}$ of $E$ is locally
finite. To simplify the notation, by definition, put: $$S_r
{:=} T_{{\xi}_r,t_r,y_r} \bigcap \left(\rk {\times} J
{\times}\rk\right) \mbox{ and }  v_r {:=}
v_{{\xi}_r,t_r,y_r} \; \; \; \; \;   \mbox{ for every } r
{\in} \nn \n{r42_44_a}$$

Furthermore, for each $\theta(\cdot) {\in} C(\rk
{\times}\rk; J),$
 and each $A {\subset} \rk {\times}\rk,$ we put by definition:
$$\Upsilon_{\theta(\cdot), A_{\theta}} = \{ (\eta,s,z) \in
\rk \times\rr\times\rk| \; s \leq \theta(\eta,z) \}
\setminus$$ $$  \{ (\eta,s,z) \in \rk \times\rr\times\rk |
\; (s= \theta(\eta,z)) \wedge((\eta,z) \in A_{\theta}) \}$$
$$ {\Gamma}_{\theta (\cdot), A_{\theta}} = \{ (\eta,s,z)
\in \rk \times\rr\times\rk| \; s \geq \theta(\eta,z) \}
\setminus$$ $$\{ (\eta,s,z) \in \rk \times\rr\times\rk | \;
(s= \theta(\eta,z)) \wedge((\eta,z) \in A_{\theta} ) \}.$$
Let ${\gothS}$ be the system of all the sets given by
$$\Sigma_{\Theta(\cdot),\vartheta(\cdot),A_{\Theta},A_{\vartheta}}
:= \Upsilon_{\Theta(\cdot), A_{\Theta}} \bigcap
{\Gamma}_{\vartheta(\cdot), A_{\vartheta}}, $$ where
$\Theta (\cdot),$ and  $\vartheta(\cdot)$ run through the
set of all the functions of class
 $C(\rk {\times}\rk; J)$ such that, for all
$(\xi, y, z) {\in} \rk{\times}\rk{\times}\rk,$
 $$ | \Theta(\xi,y) {-}\Theta(\xi,z)| {\leq} L(\xi) |y{-}z| \;  \mbox{ and } \;   |\vartheta(\xi,y) {-} \vartheta(\xi,z) | {\leq} L(\xi)  |y{-}z|,
 \;    $$
 and $A_{\Theta} {\subset} \rk{\times}\rk,$  $A_{\vartheta} {\subset} \rk{\times}\rk$
run through the set of all subsets of $\rk{\times}\rk.$

First, note that, if ${\tau}_j (\cdot,\cdot),$ $j {=}
1,...,M,$ are some functions of $\rk {\times} \rk$  to
$\rr$ such that $$\forall \xi \in \rk \;\;  \forall y \in
\rk \;\;  \forall z \in \rk \; \; \; \; \; \;\; \; \;    |
{\tau}_j (\xi,y) - {\tau}_j (\xi,z)| \leq
 L(\xi)\;  |y-z|, \; \;j = 1,...,M,    \n{r42_45}$$
then, we obtain: $$\forall \xi \in \rk \;\;  \forall y \in
\rk \;\;  \forall z \in \rk \;\; \; \;  |
\max\limits_{j=1,...,M} \{ {\tau}_j (\xi,y) \} -
\max\limits_{j=1,...,M} \{ {\tau}_j (\xi,z) \}| \leq
 L(\xi)\;  |y-z|,   \n{r42_46}$$
$$\forall \xi \in \rk \;\;  \forall y \in \rk \;\;  \forall
z \in \rk \;\; \; \;  | \min\limits_{j=1,...,M} \{ {\tau}_j
(\xi,y) \} -  \min\limits_{j=1,...,M} \{ {\tau}_j (\xi,z)
\}| \leq
 L(\xi) \; |y-z|.   \n{r42_47}$$
Therefore, it is easy to verify that
 ${\gothS}$ is a semiring of sets, i.e.,
 (see. \cite[p. 51]{kolmogorov}) first,
$\emptyset {\in} {\gothS},$ second, for each ${\Sigma}'
{\in} {\gothS},$ and each ${\Sigma}'' {\in} {\gothS},$
 we have ${\Sigma}'{\bigcap}{\Sigma}'' \in  {\gothS},$ and, third,
 for every $\Sigma {\in} {\gothS}, $ and every
$\Sigma_1 {\in} {\gothS},$ if $\Sigma_1 {\subset} \Sigma,$
then there exists a finite sequence $\{ \Sigma_i
{\}}_{i=1}^l \subset {\gothS} $ of sets from ${\gothS}$
such that $\Sigma = \bigcup\limits_{j=1}^{l} {\Sigma}_j,$
and ${\Sigma}_i {\bigcap} {\Sigma}_j = \emptyset$ for all
$i{\not=}j,$ $\{i,j\} {\subset} \{1,...,l\}.$

On the other hand, from (\ref{r42_39})-(\ref{r42_41}), we
obtain that, for each $r {\in} \nn,$ the set $S_r$ can be
represented as $S_r =  \Sigma_{{\overline \Theta}_r
(\cdot), {\overline \tau}_r (\cdot),A_{{\overline
\Theta}_r}, A_{{\overline \vartheta}_r}},$ where $$
{\overline \Theta}_r  (\eta,z) = \min \left\{ T,  \; \;
\max \{ {\theta}_{\xi_r,t_r,y_r} (\eta,z), \; t_1 \}
\right\}, \; \; \; (\eta,z) \in \rk \times \rk; $$ $$
{\overline \tau}_r  (\eta,z) = \max \left\{ t_1, \;  \;
\min \{ {\tau}_{\xi_r,t_r,y_r} (\eta,z), \; T \} \right\},
\; \; \; (\eta,z) \in \rk \times\rk; $$ $$ A_{{\overline
\Theta}_r} = \{ (\eta,z) \in \rk \times\rk | \;
{\theta}_{\xi_r, t_r,y_r} (\eta,z) \leq T \}; \; \; \;
 A_{{\overline \tau}_r} = \{ (\eta,z) \in \rk \times\rk | \; t_1 \leq {\tau}_{\xi_r, t_r,y_r} (\eta,z)  \}.$$
In addition, from (\ref{r42_39}), and from (\ref{r42_40}),
we get ${\overline \Theta}_r  (\eta,z) {=} t_1,$ and
${\overline \tau}_r (\eta,z) {=} T$ for $|\eta - \xi_r|
\geq \frac{1}{2},$ and for every $ z {\in} \rk.$ Hence,
since (\ref{r42_45}) implies (\ref{r42_46}) and
(\ref{r42_47}), from (\ref{r42_36}), and from
(\ref{r42_36_a}), it follows that, for each $ \eta {\in}
\rk,$ each $y {\in} \rk, $ and each $ z {\in} \rk,$ we
have: $$ |{\overline \Theta}_r  (\eta,y) - {\overline
\Theta}_r  (\eta,z)| \leq L(\eta)\;  |y-z|, \; \mbox{ and }
\; |{\overline \tau}_r (\eta,y)- {\overline \tau}_r
(\eta,z)| \leq L(\eta)\;  |y-z|. $$ Thus, each set $ S_r $
($r {\in} \nn$) is an element of semiring ${\gothS}.$

Then, from (\ref{r42_44}),  from  (\ref{r42_44_a}), and
from lemma 2 in \cite[p. 53]{kolmogorov}, it follows that
there exist a sequence $\{ \Sigma_l {\}}_{l=1}^{\infty} =
\{  {\Sigma}_{{\Theta}_l (\cdot), {\vartheta}_l (\cdot) ,
A_{\Theta_l} , A_{\vartheta_l}}  {\}}_{l=1}^{\infty}$
 of sets from ${\gothS}$ and a strictly
increasing sequence $\{l_q {\}}_{q=1}^{\infty} \subset \nn$
such that: ${\rm (A_1)}$ for each $q {\in} \nn$ we have
$\bigcup\limits_{r=1}^{r_q} S_r  {=}
\bigcup\limits_{l=1}^{l_q } \Sigma_l $ (which implies
$\bigcup\limits_{l=1}^{\infty} \Sigma_l =
\bigcup\limits_{r=1}^{\infty} S_r$); ${\rm (A_2)}$
${\Sigma}_{l'} \bigcap {\Sigma}_{l''} {=} \emptyset $ for
all $ l' {\not=} l'',$   $l'{\in}\nn,$  $l''{\in}\nn;$
${\rm (A_3)}$
  for each $r {\in} \nn,$ there exists
a finite set of indices $P(r) {\subset} \nn  $  such that
$S_r = \bigcup\limits_{l \in P(r)} \Sigma_l.$ Then, using
${\rm (A_1)}$ and ${\rm (A_2)}$, for each $l {\in} \nn$ we
obtain that $\Sigma_l
\subset\bigcup\limits_{r{=}r_q{+}1}^{r_{q{+}1}} S_r$
whenever $l_q{+}1 {\leq} l {\leq} l_{q{+}1},$ $q {\in}
\nn,$ and $\Sigma_l \subset \bigcup\limits_{r=1}^{r_1} S_r$
whenever $1 {\leq} l {\leq} l_1.$ Therefore, for each $l
{\in} \nn,$ there exists $r(l) {\in} \nn$ such that
$\Sigma_l {\subset} S_{r(l)}, $ and, if $1{\leq} l {\leq}
l_1,$ then $1 {\leq} r(l) {\leq} r_1, $ and if $l_q {+}1
{\leq} l {\leq} l_{q{+}1}$ $(q {\in} \nn),$ then $r_q {+}1
\leq r(l) \leq r_{q{+}1}.$ Taking into account
(\ref{r42_42}) and the inclusion $\Sigma_l \subset
T_{{\xi}_{r(l)}, t_{r(l)}, y_{r(l)}},$ we obtain: $$\left(
B_{\frac{1}{2}} (\xi) {\times} J {\times} \rk \right)
\bigcap {\Sigma}_l  = \emptyset \; \; \;  \mbox{ whenever }
\xi {\in} {\Xi}_{q{+}1}, \;\;  l {\notin} \{ l' {\}}_{l' =
l_{q-1} +1}^{l_{q+2}}, \; \;  q {\in} \nn {\setminus} \{
1\}, \; \; l {\in} \nn;  \n{r42_54}$$ $$\left(
B_{\frac{1}{2}} (\xi) {\times} J {\times} \rk \right)
\bigcap {\Sigma}_l  = \emptyset \; \; \; \mbox{ whenever }
\;  \xi \in {\Xi}_1 \bigcup {\Xi}_2 , \;\;  l \notin \{ l'
{\}}_{l' = 1}^{l_{3}}, \; \; l \in \nn.  \n{r42_55}$$
 For each $\xi {\in} \rk,$ let
$\Omega(\xi)$ be the finite number of indices given by $$
\Omega(\xi) = \left\{ \begin{array}{l} \{  l
{\}}_{l=l_{q-1}+1}^{l_{q+2}}, \; \; \; \mbox{ if } \;  \xi
\in {\Xi}_{q+1}, \;  q \in \nn,\;  q \geq 2\\ \{  l
{\}}_{l=1}^{l_{3}}, \; \; \; \; \; \; \; \; \; \mbox{ if }
\;  \xi \in {\Xi}_1 \bigcup {\Xi}_2.
\end{array}
\right. \n{r42_55_c}$$ By definition, put: $$v(\xi,t,y) =
v_{r(l)} \;\quad\quad\quad \; \; \; \; \; \; \; \mbox{
whenever } \;\;   (\xi,t,y) {\in} {\Sigma}_l, \;\;  l {\in}
\nn.  \n{r42_55_a}$$
 Then, from (\ref{r42_43}), from (\ref{r42_55_a}),
and from the inclusion $\Sigma_l \subset T_{ {\xi}_{r(l)},
t_{r(l)}, y_{r(l)}},$
 we obtain:
$$\forall (\eta,s,z) \in \bigcup\limits_{l=1}^{\infty}
{\Sigma}_l \;\; \;  \; \; \; | \dot x_p (\eta,s) -
f_p(s,z,v(\eta,s,z)) | < \delta (\eta).  \n{r42_55_b}$$ Let
us prove the following lemma.

{\bf Lemma 3.1.1.}      {\sl     1)  For each $\xi {\in}
\rk,$ there exist a unique   $z(\xi,\cdot) {\in} C(J; \rk)$
such that $$ z(\xi,T) = \xi, \n{r42_56}$$ a unique finite
sequence of indices $\{ {\nu}_j (\xi)  {\}}_{j {=}
1}^{N(\xi)} {=} \{ {\nu}_j {\}}_{j=1}^{N(\xi)} {\subset}
\Omega(\xi)$ such that $N(\xi) \leq | \Omega(\xi)|$ and
${\nu}_{\mu} {\not=} {\nu}_j$ whenever $\mu {\not=} j,$ and
a unique finite sequence $T {=} {\tau}_1^{\ast}(\xi) {>}
{\tau}_2^{\ast} (\xi) {>}...{>}{\tau}_{N(\xi)}^{\ast} (\xi)
{>}  {\tau}_{N(\xi){+}1}^{\ast} (\xi) {=}$ ${=} t_1$
 such that:

1.1)  $\dot z (\xi,t)$ is defined and continuous at each $t
{\in} J {\setminus} \{ {\tau}_1^{\ast}
(\xi),...,{\tau}_{N(\xi)}^{\ast} (\xi) \},$
 and
$$(\xi,t,z(\xi,t)) {\in} E \;  \mbox{ and } \; |\dot x_p
(\xi,t) {-} f_p (t,z(\xi,t),v(\xi,t,z(\xi,t))) | {<}
\delta(\xi),    \; \;  \mbox{ whenever } t {\in} J;
\n{r42_57}$$

1.2) for each $j{=} 1,...,N(\xi),$ we have: $$
(\xi,t,z(\xi,t)) \in {\Sigma}_{\nu_j} \; \;\;\; \; \; \; \;
\;\;\; \; \; \; \; \; \; \; \; \;  \; \; \; \; \; \; \; \;
\;  \;\;\; \; \; \;  \; \; \;  \mbox{ for all } \; \;  t
\in ]{\tau}_{j+1}^{\ast}  (\xi), {\tau}_j^{\ast} (\xi) [
,  \n{r42_58}$$ $$ \dot z (\xi,t) = \varphi
(t,z(\xi,t),v(\xi,t,z(\xi,t)))   \; \;\;\; \; \; \;\;
\;\;\; \; \; \;\; \; \mbox{ for all }  \; \;  t \in
]{\tau}_{j+1}^{\ast} (\xi), {\tau}_j^{\ast} (\xi)[  ,
\n{r42_59}$$ $${\tau}_j^{\ast} (\xi) = \Theta_{\nu_j}
(\xi,z(\xi,{\tau}_j^{\ast} (\xi))), \; \; \;\; \;\;\; \; \;
\;\; \;\;\; \; \; \; \; \; \; \;  \; {\tau}_{j+1}^{\ast}
(\xi) = {\vartheta}_{\nu_j} (\xi,z(\xi,{\tau}_{j+1}^{\ast}
(\xi))).  \n{r42_59_a}$$

2) For each $ \xi {\in} \rk,$ and each  $l {\in} \nn,$ let
   $t {\mapsto} s_l(\xi,t)$ and   $t {\mapsto} t_l (\xi,t) $
  be  given by
$$   s_l (\xi,t) {=} t {-} \vartheta_l (\xi,z(\xi,t)), \;
\; \; \;   t_l(\xi,t) {=} t {-} \Theta_l (\xi,z(\xi,t)),
\; \;  \; \; \mbox{ whenever }  t {\in} J, \;   l {\in} \nn
\n{r42_60}$$ Then, for every $ \xi {\in} \rk,$ and every
$l {\in} \nn,$ first, $$ \frac{3(t{-}\tau)}{2} \geq
s_l(\xi,t) - s_l(\xi,\tau) \geq \frac{t{-}\tau}{2} \; \;
\mbox{ whenever } \;    t{>} \tau, \; l {\in} \nn,
\n{r42_61}$$ $$ \frac{3(t{-}\tau)}{2} \geq  t_l(\xi,t) -
t_l(\xi,\tau) \geq \frac{t{-}\tau}{2}   \; \; \mbox{
whenever } \;    t{>} \tau, \; l {\in} \nn, \n{r42_62}$$
for all $t{\in} J$ and $\tau {\in} J,$ and, second, there
exist unique $s_l^{\ast} (\xi){\in}J$ and $t_l^{\ast} (\xi)
{\in}J$ such that $s_l(\xi,s_l^{\ast} (\xi)) {=} 0$ and
$t_l(\xi,t_l^{\ast} (\xi)) {=} 0.$  Moreover, $T{=}
t_{\nu_1}^{\ast} (\xi); $ ${\tau}_i^{\ast}(\xi) {=}
t_{\nu_i}^{\ast}(\xi) {=} s_{\nu_{i-1}}^{\ast} (\xi),$ $i
{=} 2,...,N(\xi);$ $t_1 {=} s_{\nu_{N(\xi)}}^{\ast}(\xi).$
}

{\bf  Proof of lemma 3.1.1.} Choose and fix an arbitrary
$\xi {\in} \rk.$ Choose  $q {\in} \nn$ such that $\xi {\in}
{\Xi}_q.$ By definition, put: ${\tau}_0^{\ast} (\xi) {=}
T,$ ${\tau}_1^{\ast} (\xi) {=}T.$ Next, using the induction
over $i {\in} \nn,$ we construct the desired
${\tau}_i^{\ast} (\xi)$ and $\nu_i {=} \nu_i (\xi),$ and
the trajectory $t {\mapsto} z(\xi,t)$ on $[ \tau_i^{\ast}
(\xi), T]$ and prove the uniqueness of the construction.
Throughout the proof of lemma 3.1.1,
 $\xi$ is assumed to be fixed, and, therefore,
we always write $\nu_j $ instead of $\nu_j (\xi).$

For $i {=} 1,$ the construction is trivial: put $z(\xi,t)
{:=} \xi$ for $t {\in} [T,T] {=} [{\tau}_1^{\ast} (\xi),
{\tau}_0^{\ast} (\xi)].$ Then, by the definition of $K_q$
and $E_q$ (see (\ref{r42_30}),
(\ref{r42_33}),(\ref{r42_34})), from the equality $y(\xi,T)
{=} \xi,$ we get: $(\xi,T,\xi) {=} (\xi,T,z(\xi,T)) {\in}
E.$ Hence, from (\ref{r42_55_b}), we get (\ref{r42_57}) for
$i{=}1$ and for $t {\in}
[{\tau}_1^{\ast}(\xi),\tau_0^{\ast} (\xi)] {=} [T,T].$
Furthermore, for $i{=}1,$ and for each $j {\in} \nn$ such
that $1{\leq} j {\leq} i{-}1,$ conditions
(\ref{r42_58})-(\ref{r42_59_a}) hold by definition (there
are $p{-}1{=}0$ identities or inclusions to be satisfied).
Finally, the uniqueness of $t {\mapsto} z(\xi,t)$ defined
on $[\tau_0^{\ast}(\xi) ,\tau_1^{\ast} (\xi)] = \{ T\}$
follows from (\ref{r42_56}), and the uniqueness of
$\tau_0^{\ast} (\xi)$ and $\tau_1^{\ast} (\xi)$ follows
from their definition.

Assume that, for some $i {\in} \nn,$ we have already
constructed a finite sequence of indices $\{ \nu_j
{\}}_{j=1}^{i-1}$ (if   $i{=}1,$ the sequence is empty by
the induction hypothesis -- see above) such that $\nu_{\mu}
{\not=} \nu_j$ for all $j {\not=} \mu$ and $\{j, \mu\}
\subset \{ 1,...,i-1\},$ a finite sequence $T{=}
\tau_0^{\ast}(\xi){=}\tau_1^{\ast}(\xi){>}\tau_2^{\ast}(\xi){>}...{>}\tau_i^{\ast}(\xi)
{\geq} t_1,$
 and a trajectory
$t {\mapsto} z(\xi,t)$ defined on $[\tau_i^{\ast}(\xi),T]$
such that:
 (\ref{r42_56}) holds,
 (\ref{r42_57}) holds whenever $t {\in} [\tau_i^{\ast}(\xi),T]$,
 and
(\ref{r42_58})-(\ref{r42_59_a}) hold whenever $1 {\leq} j
{\leq} i{-}1, $  $j {\in}\nn$ (again, for $i{=}1,$ we deal
with the empty set of conditions
(\ref{r42_58})-(\ref{r42_59_a}) -- see the induction
hypothesis). In addition, assume that we have proved the
uniqueness of $\{  \nu_j {\}}_{j=1}^{i-1}$ and $\{
\tau_j^{\ast}(\xi) {\}}_{j=0}^{i},$
 and the uniqueness of such $z(\xi,\cdot)$
 on $[\tau_i^{\ast} (\xi),T].$
Finally, suppose  we have proved that functions $s_l
(\xi,\cdot)$ and $t_l(\xi,\cdot)$ defined by (\ref{r42_60})
for all $t {\in} [\tau_i^{\ast}(\xi),T]$ and $l {\in} \nn$
satisfy
  (\ref{r42_61}), (\ref{r42_62})
for all $t \in [\tau_i^{\ast}(\xi),T]$ and $\tau \in
[\tau_i^{\ast}(\xi),T]$ (again, for $i{=}1,$ this is
trivial: in this case,  we get $[\tau_i^{\ast} (\xi),T] {=}
[T,T]$ so that the set of $\{\tau ,t \} \subset [T,T]$ such
that $t {>} \tau$ is empty). If $\tau_i^{\ast} (\xi) {=}
t_1,$ then, we put $i = : N(\xi) +1,$  and note  that, in
this case, lemma 3.1.1 is proved because, from the
definition of
 ${\gothS},$ and from
(\ref{r42_60}), we obtain: $$s_l(\xi,T) \geq 0, \; \;
t_l(\xi,T) \geq 0,\; \;  s_l(\xi,t_1) \leq 0, \; \;
t_l(\xi,t_1) \leq 0,  \; \; \; \mbox{ whenever } \; l {\in}
\nn,$$ and, therefore, from (\ref{r42_61}), (\ref{r42_62}),
we get the existence and uniqueness of $s_l^{\ast}(\xi)
{\in} J$ and $t_l^{\ast}(\xi) {\in} J$ such that
$s_l(\xi,s_l^{\ast} (\xi))=0,$ $t_l(\xi,t_l^{\ast}
(\xi))=0$ ($l \in \nn$). From the uniqueness, and from
(\ref{r42_59_a}), we get: $$T = t_{\nu_1}^{\ast}(\xi);\; \;
\;  {\tau}_i^{\ast} (\xi) = t_{\nu_i}^{\ast} (\xi) =
s_{\nu_{i-1}}^{\ast} (\xi), \; \; i = 2,...,N(\xi); \; \;
\;  t_1 = s_{\nu_{N(\xi)}}^{\ast}(\xi).$$

Therefore, it suffices to consider the case when
$\tau_i^{\ast} (\xi) {>} t_1.$ From the induction
hypothesis, we get $(\xi, \tau_i^{\ast} (\xi), z(\xi,
\tau_i^{\ast} (\xi))) {\in} E;$ hence, there exists
$\overline{\varepsilon}{>}0$ such that $(\xi, \tau_i^{\ast}
(\xi) {-} s, z(\xi, \tau_i^{\ast} (\xi))) {\in} E$
 for all $s {\in} ]0, \overline{\varepsilon}].$
Since all  $\Sigma_l$ are mutually disjoint, and $E \subset
\bigcup\limits_{l=1}^{\infty} \Sigma_l,$ we obtain that
there exist unique $\nu_i {=} \nu_{i} (\xi) {\in} \Omega
(\xi)$ and $\overline{\tau} {\in} [t_1, \tau_i^{\ast}
(\xi)[$ such that $\{ \xi \} {\times} ]\overline{\tau},
\tau_i^{\ast} (\xi)[ {\times} \{ z(\xi,\tau_i^{\ast} (\xi))
\} \subset \Sigma_{\nu_i},$ and ${\overline \tau} =
\vartheta_{\nu_i} (\xi,z(\xi,\tau_i^{\ast} (\xi))),$ ${
\tau}_i^{\ast} (\xi) =  \Theta_{\nu_i}
(\xi,z(\xi,\tau_i^{\ast} (\xi))).$
Consider the trajectory $t {\mapsto} y(t,\tau_i^{\ast}
(\xi), z(\xi, \tau_i^{\ast} (\xi)),v_{r(\nu_i)})$ (which is
defined on the maximal possible interval of time according
to our notation). It is one or the other: either it is
contained in
 $\overline{B_{d_q} (0)}$ for all  $t$  from its domain, and, then,
 it is well defined for all  $t {\in} [t_1, \tau_i^{\ast} (\xi)],$
or there exists $s {\in} ]t_1, \tau_i^{\ast} (\xi)   ]$
such that this trajectory is defined for all $ t {\in} [s,
\tau_i^{\ast} (\xi)],$ and $y (s,\tau_i^{\ast} (\xi),
z(\xi,\tau_i^{\ast} (\xi)),  v_{r(\nu_i)} ) {\not\in}
\overline{B_{d_q}(0)}.$ In the first case, by definition,
put
 $\overline{ s} {:=} t_1,$ and, in the second one, put:
$$\overline{s} := \sup \{ t \in [t_1, {\tau}_i^{\ast}
(\xi)[ | \; \; \;|y (t,\tau_i^{\ast} (\xi),
z(\xi,\tau_i^{\ast} (\xi)), v_{r(\nu_i)} ) | = d_q \} $$
(since
 $y (\cdot,\tau_i^{\ast} (\xi), z(\xi,\tau_i^{\ast} (\xi)), v_{r(\nu_i)} ) $
is continuous, we deal with the supremum of a compact
subset of $\rr,$ and
  $|y (\overline{s},\tau_i^{\ast} (\xi), z(\xi,\tau_i^{\ast} (\xi)), v_{r(\nu_i)} )| = d_q $).
In both cases, $\overline{s}$ satisfies the following
conditions: $$ \forall t \in [\overline{s}, \tau_i^{\ast}
(\xi)] \; \; \;   y (t,\tau_i^{\ast} (\xi),
z(\xi,\tau_i^{\ast} (\xi)), v_{r(\nu_i)} ) \in
\overline{B_{d_q} (0)} \n{r42_64}  $$ $$(\xi,\overline{s},
y (\overline{s},\tau_i^{\ast} (\xi), z(\xi,\tau_i^{\ast}
(\xi)), v_{r(\nu_i)} )) \notin {\rm int} \Sigma_{\nu_i}.
\n{r42_65}$$ Let us prove that the functions ${\overline
s}_l (\cdot,\cdot)$ and ${\overline t}_l (\cdot,\cdot),$
$l{\in}\nn,$ defined by $${\overline s}_l (\xi,t) {=} t {-}
{\vartheta}_l (\xi, y (t,\tau_i^{\ast} (\xi),
z(\xi,\tau_i^{\ast} (\xi)), v_{r(\nu_i)} )),
 \; {\overline t}_l (\xi,t) {=} t {-} {\Theta}_l (\xi, y (t,\tau_i^{\ast} (\xi), z(\xi,\tau_i^{\ast} (\xi)), v_{r(\nu_i)} ))$$
satisfy the conditions: $$ \forall  \tau \in [\overline{s},
\tau_i^{\ast} (\xi) ] \; \;  \;  \forall t \in ] \tau,
\tau_i^{\ast} (\xi)] \; \; \;  \forall l \in \nn\; \; \;
\; \; \;
 \left( \frac{3(t-\tau)}{2} \geq {\overline s}_l (\xi,t) - {\overline s}_l (\xi,\tau )  \geq \frac{t-\tau}{2} \right) \wedge$$
$$ \wedge \left(\frac{3(t-\tau)}{2} \geq {\overline t}_l
(\xi,t) - {\overline t}_l (\xi,\tau ) \geq
\frac{t-\tau}{2}\right). \n{r42_65_a} $$ Indeed, since all
$\Sigma_l$ are elements of ${\gothS},$ from (\ref{r42_64}),
(\ref{r42_36}),  (\ref{r42_36_a}),
 and from the definition of ${\gothS},$
it follows that, for every $l {\in} \nn,$ every $\tau {\in}
[\overline{s}, \tau_i^{\ast} (\xi)],$ and every $t \in
]\tau, \tau_i^{\ast} (\xi)],$ we obtain $$ | \vartheta_l
\left( \xi, y (t,\tau_i^{\ast} (\xi), z(\xi,\tau_i^{\ast}
(\xi)), v_{r(\nu_i)} ) \right) - \vartheta_l \left( \xi, y
(\tau,\tau_i^{\ast} (\xi), z(\xi,\tau_i^{\ast} (\xi)),
v_{r(\nu_i)} ) \right) | \leq$$ $$\leq L(\xi)|  y \left(
t,\tau_i^{\ast} (\xi), z(\xi,\tau_i^{\ast} (\xi)),
v_{r(\nu_i)} \right) - y \left( \tau,\tau_i^{\ast} (\xi),
z(\xi,\tau_i^{\ast} (\xi)), v_{r(\nu_i)} \right)  |\leq$$
$$\leq L(\xi) \int\limits_{\tau}^t |\varphi \left( s, y
(s,\tau_i^{\ast} (\xi), z(\xi,\tau_i^{\ast} (\xi)),
v_{r(\nu_i)} ), v_{r(\nu_i)} \right) | ds <
\frac{t-\tau}{2}.$$ From this, we get the first group of
inequalities
 (\ref{r42_65_a}) for functions
  ${\overline s}_l (\cdot,\cdot).$ The proof of the the inequalities (\ref{r42_65_a})
for ${\overline t}_l (\cdot,\cdot)$ is similar.

Note that, by the construction, ${\overline t}_{\nu_i}
(\xi, \tau_i^{\ast} (\xi)) {=} 0,$ and ${\overline
s}_{\nu_i} (\xi, \tau_i^{\ast} (\xi)){>}0.$ From
(\ref{r42_65_a}), we get, in particular, that ${\overline
t}_{\nu_i} (\xi, \cdot)$ and ${\overline s}_{\nu_i} (\xi,
\cdot)$ are strictly increasing functions on
$[\overline{s}, \tau_i^{\ast} (\xi)]; $ hence ${\overline
t}_{\nu_i} (\xi, \overline{s}) {<}0.$ Then, the inequality
${\overline s}_{\nu_i} (\xi, \overline{s}) {>}0$ is
impossible, because it implies $$\vartheta_{\nu_i} \left(
\xi, y (\overline{s},\tau_i^{\ast} (\xi),
z(\xi,\tau_i^{\ast} (\xi)), v_{r(\nu_i)} ) \right) <
\overline{s} < \Theta_{\nu_i} \left( \xi, y
(\overline{s},\tau_i^{\ast} (\xi), z(\xi,\tau_i^{\ast}
(\xi)), v_{r(\nu_i)} ) \right),$$
 which contradicts (\ref{r42_65}). Therefore
${\overline s}_{\nu_i} (\xi,\overline{s}) \leq 0.$ Then,
since ${\overline s}_{\nu_i} (\xi,\cdot)$ is strictly
increasing and continuous on $[\overline{s},\tau_i^{\ast}
(\xi)],$ there exists a unique $\tau_{i{+}1}^{\ast} (\xi)
{\in} [\overline{s}, \tau_i^{\ast} (\xi)]$ such that
${\overline s}_{\nu_i} (\xi,\tau_{i+1}^{\ast} (\xi)) {=}
0,$
 and such that for each
$t {\in} ]\tau_{i+1}^{\ast} (\xi), \tau_{i}^{\ast} (\xi)[$
we have $ {\overline t}_{\nu_i} (\xi,t){<} 0{<}{\overline
s}_{\nu_i} (\xi,t),$ i.e., (by the definition of
${\overline s}_{l} (\xi,\cdot)$ and ${\overline t}_{l}
(\xi,\cdot)$) $$ \forall t  \in ]{\tau}_{i+1}^{\ast}
(\xi),{\tau}_i^{\ast} (\xi)[  \; \; \; \; \; \; \left(
\xi,t,y (t,\tau_i^{\ast} (\xi), z(\xi,\tau_i^{\ast} (\xi)),
v_{r(\nu_i)} ) \right) \in {\rm int} \Sigma_{\nu_i};
\n{r42_67}  $$ $$\tau_{i+1}^{\ast} (\xi) =
\vartheta_{\nu_i}\left(\xi,y
(\tau_{i+1}^{\ast}(\xi),\tau_i^{\ast} (\xi),
z(\xi,\tau_i^{\ast} (\xi)), v_{r(\nu_i)} ) \right), $$ $$
\tau_i^{\ast}(\xi) =\Theta_{\nu_i} \left(\xi,y
(\tau_{i}^{\ast}(\xi),\tau_i^{\ast} (\xi),
z(\xi,\tau_i^{\ast} (\xi)), v_{r(\nu_i)} )\right).
\n{r42_68}$$ Define the extension of $z(\xi,\cdot)$ to
$[t_{i{+}1}^{\ast} (\xi), T]$ by $$z(\xi,t) = y
(t,\tau_i^{\ast} (\xi), z(\xi,\tau_i^{\ast} (\xi)),
v_{r(\nu_i)} )\; \; \; \;  \mbox{ whenever }\; \; t \in
[\tau_{i+1}^{\ast} (\xi),\tau_{i}^{\ast} (\xi)].
\n{r42_69}$$ Then, from (\ref{r42_67}) and (\ref{r42_68}),
we obtain that (\ref{r42_58})-(\ref{r42_59_a}) hold not
only for $1{\leq}j{\leq}i{-}1$ but for all $j{=}1,...,i.$
Furthermore, from (\ref{r42_67}), (\ref{r42_69}),
(\ref{r42_55_b}), and from the condition $$|\dot x_p
(\xi,t) - f_p(t,z(\xi,t),v(\xi,t,z(\xi,t)))| < \delta(\xi)
\; \mbox{ whenever } t \in [{\tau}_i^{\ast} (\xi),T],$$
(which holds by the induction hypothesis) we get $$|\dot
x_p (\xi,t) - f_p(t,z(\xi,t),v(\xi,t,z(\xi,t)))| <
\delta(\xi) \;  \mbox{ whenever } t \in
[{\tau}_{i+1}^{\ast} (\xi),T].$$ By the construction, $t
{\mapsto} z(\xi,t)$ is the trajectory, of (\ref{r4_25}),
defined on $[{\tau}_{i+1}^{\ast} (\xi),T]$ by the initial
condition $z(\xi,T) {=} \xi$ and by the control $t
{\mapsto} f_p (t,z(\xi,t),v(\xi,t,z(\xi,t)));$ therefore,
from lemma 3.3, we obtain that
  $z(\xi,t) {\in} K_q$ for all $t {\in} [{\tau}_{i+1}^{\ast} (\xi),T];$
hence, by the induction hypothesis, $(\xi,t,z(\xi,t)) \in E
$ whenever $  t \in [{\tau}_{i+1}^{\ast} (\xi),T]. $

Thus, the trajectory $z(\xi,\cdot)$ (which is now defined
on $[{\tau}_{i{+}1}^{\ast} (\xi),T]$)  satisfies
(\ref{r42_57}) for all $t {\in} [{\tau}_{i{+}1}^{\ast}
(\xi),T].$ In addition, by the induction hypothesis,
(\ref{r42_56}) holds, and, from (\ref{r42_65_a}) and
(\ref{r42_69}), it follows that conditions (\ref{r42_61}),
(\ref{r42_62}) hold for all $t {\in} [{\tau}_{i+1}^{\ast}
(\xi),T]$  and all $\tau {\in} [{\tau}_{i+1}^{\ast}
(\xi),T].$ Hence, from (\ref{r42_59_a}), we get $\nu_i
{\not=} \nu_j$ whenever $j{=}1,...,i{-}1.$ (Because, for
each $j{=}1,...,i{-}1,$ and  each $t {\in}
]{\tau}_{i{+}1}^{\ast} (\xi),\tau_i^{\ast}(\xi)[,$ we have:
$s_{\nu_j} (\xi,t) {<}0 ,$  i.e., $(\xi,t,z(\xi,t)) \not\in
\Sigma_{\nu_j},$ but $ (\xi,t,z(\xi,t)) \in
\Sigma_{\nu_i}$).

Let us prove that, $\nu_i {=} \nu_i(\xi),$
$\tau_{i{+}1}^{\ast} (\xi),$ and the trajectory
$t{\mapsto}z(\xi,t)$ for $t{\in}[{\tau}_{i+1}^{\ast}
(\xi),\tau_i^{\ast}(\xi)]$ are uniquely determined by this
procedure. Since, by the induction hypothesis, $\{  \nu_j
{\}}_{j{=}1}^{i{-}1},$ $\{  \tau_j^{\ast} (\xi)
{\}}_{j=0}^i,$ and $z(\xi,t)$ for $t{\in}[{\tau}_{i}^{\ast}
(\xi),T]$ are already uniquely determined, this means the
uniqueness of $\{  \nu_j  {\}}_{j=1}^{i},$ $\{
\tau_j^{\ast} (\xi) {\}}_{j=0}^{i+1},$
 and  $z(\xi,t)$ for all $t{\in}[{\tau}_{i{+}1}^{\ast} (\xi),T].$

Assume that there exist $\overline{t} {\in}
[t_1,\tau_i^{\ast}(\xi)[,$ ${\tilde z} (\xi,\cdot) \in
C([\overline{t}, \tau_i^{\ast} (\xi)] ; \rk) {\bigcap} C^1
(]\overline{t}, \tau_i^{\ast} (\xi)[; \rk),$
 and ${\tilde \nu} {\in} \nn$ such that
$$(\xi,t,{\tilde z}(\xi,t)) \in \Sigma_{\tilde \nu} \; \;
\; \; \; \; \; \; \; \; \; \; \; \; \; \; \; \; \;  \; \;
\; \; \; \; \; \; \; \; \; \; \; \; \mbox{ whenever } \;  t
\in ]\overline{t}, \tau_i^{\ast} (\xi)[; \n{r42_70}$$
$$\frac{d}{dt} {\tilde z}(\xi,t) = \varphi (t,{\tilde
z}(\xi,t),v(\xi,t,{\tilde z}(\xi,t))) \; \; \mbox{ whenever
} \;  t \in ]\overline{ t} , \tau_i^{\ast} (\xi)[;
\n{r42_71}$$ $$ {\tilde z} (\xi,\tau_i^{\ast} (\xi)) =
z(\xi,\tau_i^{\ast} (\xi)); \;\quad\quad\quad \; \; \; \;
\; \; \; \; \; \; \; \; \; \; \; \; \; \; \; \; \; \; \; \;
\; \; \; \; \; \; \; \; \; \; \; \; \;   \n{r42_72}$$ $$
\vartheta_{\tilde{\nu}} (\xi,{\tilde z} (\xi,\overline{t}))
= \overline{t}.\;\; \;\quad\quad\quad \; \; \; \; \; \; \;
\; \; \; \; \; \; \; \; \; \; \; \; \; \; \; \; \; \; \; \;
\; \; \; \; \; \;\;  \; \;\quad\quad\quad   \n{r42_73}$$

Let us prove that, from (\ref{r42_70})-(\ref{r42_73}), it
follows that ${\tilde \nu} {=} { \nu}_i,$ $\overline{t} {=}
\tau_{i{+}1}^{\ast} (\xi),$ and ${\tilde z} (\xi,t) {=}
z(\xi,t)$ whenever $t \in [\tau_{i+1}^{\ast} (\xi),
\tau_i^{\ast}(\xi)].$ Indeed, for each $t {\in}
[\overline{t},\tau_i^{\ast}(\xi)]$ we obtain that $$
{\tilde z}(\xi,t) = z(\xi, \tau_i^{\ast} (\xi)) -
\int\limits_{t}^{{\tau}_i^{\ast} (\xi) } \varphi ( s,
{\tilde z}(\xi,s), v(\xi,s, {\tilde z}(\xi,s)))ds;$$ in
addition, from (\ref{r42_54}), (\ref{r42_55}), and from
 (\ref{r42_70}), we obtain that
${\tilde \nu} \in \Omega(\xi),$ which implies ${\tilde
z}(\xi,t) {\in} \overline{B_{d_q} (0)}$ for all $t {\in}
[\overline{t}, \tau_i^{\ast} (\xi)];$ hence, (by the same
argument as for the proof of
 (\ref{r42_65_a}))
 from (\ref{r42_70}), and from the definition of
$\{\Sigma_l {\}}_{l=1}^{\infty}$ and ${\gothS},$ we obtain
that,
 for every $ t \in [\overline{t}, \tau_i^{\ast} (\xi)]$
and every $ \tau\in ]\overline{t},t],$ $$|
\vartheta_{\nu_i} (\xi, {\tilde z} (\xi,t))
{-}\vartheta_{\nu_i} (\xi, {\tilde z} (\xi,\tau))| {<}
\frac{t{-}\tau}{2},
 \; \; \; \; | \Theta_{\nu_i} (\xi, {\tilde z} (\xi,t)) {-} \Theta_{\nu_i} (\xi, {\tilde z} (\xi,\tau))| {<} \frac{t{-}\tau}{2}.$$
Therefore, the continuous functions $t {\mapsto} {\tilde
s}_{\nu_i}(\xi,t) $ and $t {\mapsto} {\tilde
t}_{\nu_i}(\xi,t) $ given by $${\tilde s}_{\nu_i}(\xi,t) =
t- \vartheta_{\nu_i} (\xi, \tilde{z} (\xi,t)),\; \; \;
\;{\tilde t}_{\nu_i}(\xi,t) = t- \Theta_{\nu_i} (\xi,
\tilde{z} (\xi,t)) , \; \; \; \; \; \; \; \;  \; \mbox{ for
all } t \in [\overline{t},{\tau}_i^{\ast} (\xi)]$$ are
strictly increasing on $[\overline{t}, \tau_i^{\ast}
(\xi)].$ From (\ref{r42_72}), it follows that ${\tilde
s}_{\nu_i} (\xi,\tau_i^{\ast} (\xi)) {>}0,$ and ${\tilde
t}_{\nu_i} (\xi,\tau_i^{\ast} (\xi)) {=}0;$ hence there
exists $\tilde{ s} {\in} ]\overline{t}, \tau_i^{\ast}
(\xi)[$ such that, for each $t \in [\tilde{s},\tau_i^{\ast}
(\xi)[,$ we have ${\tilde t}_{\nu_i} (\xi,t) {<} 0 {<}
{\tilde s}_{\nu_i} (\xi,t), $ i.e., $(\xi,t,{\tilde
z}(\xi,t)) {\in} {\rm int} {\Sigma}_{\nu_i}.$ Since
$\tilde{z}(\xi,\cdot)$ is continuous, there exists $\tau
{\in} [{\overline t} , {\tilde s}[$ such that, for each $t
\in ]\tau,\tilde{s}],$ we have $(\xi,t,\tilde{z}(\xi,t))
\in {\rm int} \Sigma_{\nu_i}.$ By definition, put: $$
{\tau}^{\ast} = \inf \{ \tau \in [\overline{t},\tilde{s}[\;
|\; \; \;  \; \forall t \in ]\tau,\tilde{s}] \; \;
(\xi,t,\tilde{z}(\xi,t)) \in {\rm int} \Sigma_{\nu_i}\}.
$$ From the definition of $\tau^{\ast},$ it follows that,
for each $t \in ] \tau^{\ast}, \tilde{s}],$ we get the
inclusion $(\xi,t,\tilde{z}(\xi,t)) {\in} {\rm int}
\Sigma_{\nu_i}.$ This inclusion holds for all $t {\in}
[\tilde{s}, \tau_i^{\ast}(\xi)[$ as well (by the definition
of ${\tilde s}$); hence, for each $t \in ]\tau^{\ast} ,
\tau_i^{\ast} (\xi)[,$ we obtain $v(\xi,t,\tilde{z}
(\xi,t)) {=} v_{r(\nu_i)},$ i.e., taking into account
(\ref{r42_71}) and (\ref{r42_72}), we get $\tilde{z}
(\xi,t){=} y(t,\tau_i^{\ast} (\xi),
z(\xi,\tau_i^{\ast}(\xi)),v_{r(\nu_i)})$ for all $t {\in}
[\tau^{\ast},\tau_i^{\ast} (\xi)].$ From this,  from
(\ref{r42_70}), and from the fact that all $\Sigma_l,$
$l\in \nn,$ are mutually disjoint,
 we get  $\tilde{\nu} = \nu_j.$

Let us prove that $\tau^{\ast} {=} \tau_{i+1}^{\ast} (\xi)
{=} \overline{t}.$ First, we prove that
 $\tau^{\ast} {=} \tau_{i+1}^{\ast} (\xi).$
Indeed, if $\tau^{\ast} {<} \tau_{i+1}^{\ast} (\xi) ,$
then, by the above, the point $(\xi,\tau_{i+1}^{\ast}
(\xi), \tilde{z} (\xi, \tau_{i+1}^{\ast} (\xi))) =
(\xi,\tau_{i+1}^{\ast} (\xi), {z} (\xi, \tau_{i+1}^{\ast}
(\xi))) $ belongs to $\partial {\Sigma}_{\nu_i},$ which
contradicts the definition of $\tau^{\ast}.$ If $
\tau_{i+1}^{\ast} (\xi) {<} \tau^{\ast} {<} \tau_{i}^{\ast}
(\xi), $ then, since $t_{\nu_i} (\xi,\cdot)$ and $s_{\nu_i}
(\xi,\cdot)$ are strictly increasing on $ [
\tau_{i+1}^{\ast} (\xi) , \tau_{i}^{\ast} (\xi) ],$ we get
$t_{\nu_i} (\xi,\tau^{\ast}){<}0{<}s_{\nu_i}
(\xi,\tau^{\ast}),$ i.e., the point
$(\xi,\tau^{\ast},z(\xi,\tau^{\ast})){=}
(\xi,\tau^{\ast},{\tilde z} (\xi,\tau^{\ast}))$ belongs to
${\rm int} \Sigma_{\nu_i}.$ Then, from (\ref{r42_73}), from
the equality
 ${\tilde \nu} {=} \nu_i,$
and from the definition of $\tau^{\ast}$ (${\tau}^{\ast}
\in [\overline{t},\tilde{s}]$) it follows that
$\overline{t} {<} \tau^{\ast}.$ Therefore, since ${\tilde
z} (\xi,\cdot)$ is continuous, there exists
${\theta}^{\ast}  {\in} [\;  \overline{t},{\tau}^{\ast}[$
 such that, for each
$t {\in} ]\theta^{\ast}, \tau^{\ast}], $ and, moreover, for
each $t {\in} ]\theta^{\ast}, \tilde{s}],$ we get
$(\xi,t,\tilde{z}(\xi,t)) {\in} {\rm int} \Sigma_{\nu_i}.$
This contradicts the definition of
  $\tau^{\ast}.$ Thus,
$\tau^{\ast} {=} \tau_{i+1}^{\ast}(\xi),$ as desired.

Second, we prove that $\tau_{i+1}^{\ast} (\xi) {=}
{\tau}^{\ast} {=} \overline{t}.$ Assume the converse, then
from the definition of
 $\tau^{\ast},$ and from the equality
 $\tau_{i+1}^{\ast} (\xi) = \tau^{\ast}$
(which is proved) we get $\overline{t} {<}   {\tau}^{\ast}
{=}\tau_{i+1}^{\ast} (\xi). $ Therefore, having proved that
${\tilde s}_{\nu_i}(\xi,\cdot)$ and ${\tilde
t}_{\nu_i}(\xi,\cdot)$
 are strictly increasing on $[\overline{t},{\tau}_{i}^{\ast} (\xi)],$
and $\tilde{\nu} = \nu_i,$ from (\ref{r42_73}), we get:
$$\tau_{i+1}^{\ast}(\xi) > \vartheta_{\nu_i}(\xi,\tilde{z}
(\xi,\tau_{i+1}^{\ast}
(\xi)))=\vartheta_{\nu_i}(\xi,\tilde{z} (\xi,\tau^{\ast}
))= \vartheta_{\nu_i}(\xi,{z} (\xi,\tau^{\ast} ))= =
\vartheta_{\nu_i}(\xi,{z} (\xi,\tau_{i+1}^{\ast} (\xi))).$$
This contradicts the definition of $\tau_{i+1}^{\ast}
(\xi).$

Thus, it is proven that $\overline{t} {=}
{\tau}_{i+1}^{\ast} (\xi) {=}{ \tau}^{\ast},$ $\tilde{\nu}
{=} {\nu}_i,$ and $z(\xi,t) {=} \tilde{z}(\xi,t)$ whenever
$t \in [{\tau}_{i+1}^{\ast} (\xi),  {\tau}_i^{\ast}
(\xi)],$ as desired.

To this end, we construct the uniquely determined sequences
$ \{ \tau_i^{\ast}(\xi)  \},$ $\{ \nu_{i-1} (\xi) \},$ and
the trajectory  $t {\mapsto} z(\xi,t)$ on $[{\tau}_i^{\ast}
(\xi), {\tau}_{i-1}^{\ast}  (\xi)]$ by induction over
$i=1,2,3,....$ It is one or the other: either, for some
$i{=:}N(\xi),$ we get $t_{i+1}^{\ast}(\xi) {=} t_1,$ and,
then, lemma 3.1.1 is proved, or we obtain an infinite
sequence $\{ \tau_i^{\ast} (\xi) \}$ such that, for each
 $i {\in} \nn$ we have:  $\tau_i^{\ast} (\xi) >t_1.$
But the second case is impossible, because, by the
construction, each $\nu_i (\xi)$ belongs to the finite set
$\Omega (\xi),$ and
 $\nu_i (\xi) {\not=} \nu_j$  whenever   $i {\not=}j.$
Thus, for the obtained sequence of ${\tau}_i^{\ast} (\xi),$
 there exists $i {=} N(\xi)+1 {<} |\Omega (\xi)|$
such that $\tau_{i}^{\ast} (\xi) {=} t_1.$ The proof of
lemma 3.1.1 is complete.



Define the desired family of controls $\{ v(\xi,\cdot)
{\}}_{\xi \in \rk}$ by $$v(\xi,t) = v(\xi,t,z(\xi,t))
\quad\quad\quad \; \; \;  \; \; \mbox{ whenever }  \;   t
\in J,  \; \;  \xi \in \rk. \n{r42_74} $$ Let us prove that
this family satisfies conditions 1)- 3) of lemma 3.4.

From lemma 3.1.1 (see  (\ref{r42_56}), (\ref{r42_57}), and
(\ref{r42_59})) it follows that $\{ v(\xi,\cdot) {\}}_{\xi
\in \rk}$ satisfies condition 2) of lemma 3.4. Let
 $M(\cdot)$ be an arbitrary function of class $C(\rk; ]0, + \infty[)$
such that $M(\xi) \geq N_0 (q +3) +1$ for every $\xi {\in}
\rk,$ where $q {\in} \zplus$ is such that $\xi {\in}
{\Xi}_{q{+}1},$ and
 $N_0 (\overline{q})$  ($\overline{q} {\in} \nn$) is defined in (\ref{r42_35}).
  Then, from the definition of  $\Omega(\xi),$
  (see  (\ref{r42_55_c})) it follows that
we obtain $v_{r(l)} \in U_{q +3}$ for all $l {\in} \Omega
(\xi);$
 hence $|v_{r(l)}| {\leq} M(\xi),$
and, therefore,
   $\{ v(\xi,\cdot) {\}}_{\xi \in \rk}$ satisfies condition
 $3)$  of lemma 3.4. To conclude the
proof of lemma 3.4, it suffices to show that
 $\{ v(\xi,\cdot) {\}}_{\xi \in \rk}$ satisfies condition 1).

 Take any $\xi {\in} \rk,$
and let $\{\nu_i {=} \nu_i (\xi) {\}}_{i=1}^{N(\xi)}$ be
the sequence, of indices, defined by $\xi$ in lemma 3.1.1.
From now on, up to the very end of the proof of lemma 3.4,
$\xi$ is assumed to be fixed, and, therefore, we write
$\nu_i $ instead of   $ \nu_i (\xi)$. From lemma 3.1.1, it
follows that to prove that $$ \parallel v(\eta,\cdot) -
v(\xi,\cdot) {\parallel}_{L_1 (J; \rmpplusl)} \rightarrow 0
\; \;\; \; \; \;  \mbox{ as }   \eta \rightarrow \xi
\n{r42_75}$$ it suffices to prove that all the fuctions
$\eta {\mapsto} s_{\nu_i}^{\ast} (\eta)$ and $\eta
{\mapsto} t_{\nu_i}^{\ast} (\eta)$, where $i {\in} \{
1,...,N(\xi) \}$, are continuous at $\xi.$ Indeed,
combining this with the inequalities
 $ s_{\nu_i}^{\ast} (\xi) {<}  t_{\nu_i}^{\ast} (\xi)$,
we get the existence of $\overline{\delta}>0$ such that,
for each $\eta {\in} B_{\overline{\delta}} (\xi),$ and each
$i {\in} \{ 1,...,N(\xi) \},$ we have: $ s_{\nu_i}^{\ast}
(\eta) {<}  t_{\nu_i}^{\ast} (\eta).$ Then, from
(\ref{r42_61})-(\ref{r42_62}), we obtain that, for each
   $\eta {\in} B_{\overline{\delta}} (\xi),$
and each $i {\in} \{ 1,...,N(\xi) \},$ the inclusion
$(\eta,t, z(\eta,t)) {\in} {\rm int} \Sigma_{\nu_i}$ holds
if and only if
 $ s_{\nu_i}^{\ast} (\eta) {<} t{<} t_{\nu_i}^{\ast} (\eta);$
hence $$ v(\eta,t) = v_{r_(\nu_i)}  \; \; \;\; \; \; \;
\mbox{ whenever } \; \; \eta {\in} B_{\overline{\delta}}
(\xi),  \; \;   s_{\nu_i}^{\ast} (\eta) {<} t{<}
t_{\nu_i}^{\ast} (\eta),\; \; i= 1,...,N(\xi).$$ Then,
since  $ s_{\nu_i}^{\ast} (\cdot)$ and
   $t_{\nu_i}^{\ast} (\cdot)$ are continuous at $\xi,$
and
 $\{  v(\eta,\cdot) {\}}_{\eta \in \rk}$
satisfies condition 3) of lemma 3.4, we get (\ref{r42_75}).

 For $\nu_0 {=} 0 \in \zplus$,
by definiton, put: $$\Sigma_{\nu_0} = \emptyset,\; \;
\Theta_{\nu_0} (\eta,z) = \vartheta_{\nu_0} (\eta,z) = T,\;
\; A_{{\Theta}_{\nu_0}} = A_{{\vartheta}_{\nu_0}}= \rk
\times\rk,$$ $$  s_{\nu_0}^{\ast} (\eta) = t_{\nu_0}^{\ast}
(\eta)= {\tau}_0^{\ast} (\eta)= {\tau}_1^{\ast} (\eta) = T,
\; \; \; \; \; \; (\eta,z) \in \rk \times\rk.$$
 Let us prove by induction on $i {\in} \{ 0,1,...,N(\xi) \}$
that $ \eta {\mapsto} t_{\nu_i}^{\ast} (\eta), $  $\eta
{\mapsto} s_{\nu_i}^{\ast} (\eta),$
 $ \eta {\mapsto} z(\eta,t_{\nu_i}^{\ast} (\eta)), $ and
$\eta \mapsto z(\eta, s_{\nu_i}^{\ast} (\eta))$ are
continuous at $\xi {\in} \rk.$
 For $i{=}0$ the statement is trivial;
assume that it is proven for some
 $i {\in} \{0,1,...,N(\xi){-}1 \}.$

First, let us prove that $ \eta \mapsto
t_{\nu_{i+1}}^{\ast} (\eta) $ is continuous at $\xi$. By
definition, put $\overline{\vartheta} {:=} s_{\nu_{i+1}}
(\xi, \tau_{i+1}^{\ast} (\xi)),$
 i.e., (see lemma 3.1.1)
$$\overline{\vartheta} =  \tau_{i+1}^{\ast} (\xi) -
{\vartheta}_{\nu_{i+1}} \left( \xi, z(\xi,\tau_{i+1}^{\ast}
(\xi))\right) = {\Theta}_{\nu_{i+1}} \left( \xi,z(\xi,
\tau_{i+1}^{\ast} (\xi))\right)-$$
$$-{\vartheta}_{\nu_{i+1}} \left(\xi,
z(\xi,\tau_{i+1}^{\ast} (\xi)) \right)=
{\vartheta}_{\nu_{i}} \left( \xi, z(\xi,\tau_{i+1}^{\ast}
(\xi))\right)- {\vartheta}_{\nu_{i+1}} \left( \xi,
z(\xi,\tau_{i+1}^{\ast} (\xi))\right) >0 \n{r42_76}$$ Take
any $\varepsilon {\in} ]0,
\frac{\overline{\vartheta}}{2}].$ Since
$\vartheta_{\nu_{i+1}} (\cdot,\cdot),$ $\Theta_{\nu_{i+1}}
(\cdot,\cdot),$  $\vartheta_{\nu_{i}} (\cdot,\cdot),$
 and $\Theta_{\nu_{i}} (\cdot,\cdot)$
are continuous on  $\rk {\times}\rk,$ and
 $s_{\nu_i}^{\ast} (\cdot)$ is continuous at $\xi$,
 there exist ${\overline \delta}_1 >0$
 and ${\overline \delta}_2 >0$ such that
$$|\vartheta_{\nu_j} (\eta,y) - \vartheta_{\nu_j}
(\xi,z(\xi,\tau_{i+1}^{\ast} (\xi)))| <
\frac{\varepsilon}{4} \leq \frac{\overline{\vartheta}}{8}
\; \;\; \; \; \;  \mbox{ for all }$$ $$(\eta,y) \in
B_{{\overline \delta}_1} (\xi) \times B_{{\overline
\delta}_2} (z(\xi,\tau_{i+1}^{\ast} (\xi) )),  \; \; \;  j
\in \{ i, \; i+1\};   \n{r42_77}$$ $$|\Theta_{\nu_j}
(\eta,y) - \Theta_{\nu_j} (\xi,z(\xi,\tau_{i+1}^{\ast}
(\xi)))| < \frac{\varepsilon}{4} \leq
\frac{\overline{\vartheta}}{8}  \; \;\; \; \; \;  \mbox{
for all } $$ $$(\eta,y) \in B_{{\overline \delta}_1} (\xi)
\times B_{{\overline \delta}_2} (z(\xi,\tau_{i+1}^{\ast}
(\xi) )), \; \; \;  j \in \{ i, \; i+1\};  \n{r42_78}$$
$$|s_{\nu_i}^{\ast} (\eta) - s_{\nu_i}^{\ast} (\xi) | <
\frac{\varepsilon}{2},\; \; \;\; \;  {\vartheta}_{\nu_i}
(\eta,y) \leq {\Theta}_{\nu_i} (\eta,y),\; \; \; \; \;
{\Theta}_{\nu_{i+1}} (\eta,y) \leq {\Theta}_{\nu_i}
(\eta,y)$$ $$ \mbox{ for all } \;  (\eta,y) \in
B_{{\overline \delta}_1} (\xi) \times B_{{\overline
\delta}_2} (z(\xi,\tau_{i+1}^{\ast} (\xi) )).
\n{r42_79}$$ Since  $\Sigma_{\nu_i} \bigcap
\Sigma_{\nu_{i+1}} = \emptyset,$ from
 (\ref{r42_76})-(\ref{r42_79}), we get:
$$ {\vartheta}_{\nu_{i+1}} (\eta,y) < {\Theta}_{\nu_{i+1}}
(\eta,y)\leq {\vartheta}_{\nu_{i}} (\eta,y) \leq
{\Theta}_{\nu_i} (\eta,y)\; \; \; \; \; \; \;$$ $$ \mbox{
for all }   (\eta,y) \in B_{{\overline \delta}_1} (\xi)
\times B_{{\overline \delta}_2} (z(\xi,\tau_{i+1}^{\ast}
(\xi) ))  \n{r42_80}$$
 On the other hand, since
$\eta  {\mapsto} z(\eta,s_{\nu_i}^{\ast} (\eta))$ is
continuous at $\xi,$  and $\tau_{i+1}^{\ast} (\xi) {=}
s_{\nu_i}^{\ast} (\xi),$ we obtain that there exists
 $\overline{\delta} {\in} ] 0, {\overline\delta}_1]$
such that, for each
 $\eta {\in} B_{\overline \delta} (\xi),$
we have $z(\eta,s_{\nu_i}^{\ast} (\eta)) \in B_{{\overline
\delta}_2} (z(\xi, \tau_{i+1}^{\ast} (\xi))).$ From this,
and from (\ref{r42_80}), and (\ref{r42_76})-(\ref{r42_78}),
we obtain: $$ \forall \eta \in B_{{\overline \delta}} (\xi)
\; \; \; \; \; {\vartheta}_{\nu_{i+1}} (\eta, z(\eta,
s_{\nu_i}^{\ast} (\eta))) +\varepsilon <
{\vartheta}_{\nu_{i}} (\eta, z(\eta, s_{\nu_i}^{\ast}
(\eta))) - \frac{\varepsilon}{2}<$$ $$<
{\Theta}_{\nu_{i+1}} (\eta, z(\eta, s_{\nu_i}^{\ast}
(\eta))) \leq {\vartheta}_{\nu_{i}} (\eta, z(\eta,
s_{\nu_i}^{\ast} (\eta))) = s_{\nu_i}^{\ast}(\eta).$$
Combining this with (\ref{r42_61}), (\ref{r42_62}), we
obtain: $$ s_{\nu_i}^{\ast} (\eta) - \varepsilon -
{\Theta}_{\nu_{i+1}} (\eta, z(\eta, s_{\nu_i}^{\ast} (\eta)
{-} \varepsilon)) <0, \quad\; \;  s_{\nu_i}^{\ast} (\eta) -
\varepsilon - {\vartheta}_{\nu_{i+1}} (\eta, z(\eta,
s_{\nu_i}^{\ast} (\eta) {-} \varepsilon)) >0,$$ i.e.,
$\left( \eta, s_{\nu_i}^{\ast} (\eta) {-} \varepsilon,
z(\eta, s_{\nu_i}^{\ast} (\eta) {-} \varepsilon) \right)
{\in} {\rm int} \Sigma_{\nu_{i+1}}.$
 By lemma 3.1.1, functions $s_{\nu_{i+1}} (\eta,\cdot)$ and
$t_{\nu_{i+1}} (\eta,\cdot)$ are strictly increasing;
hence, for each $\eta {\in} B_{{\overline \delta}} (\xi),$
we have $t_{\nu_{i+1}}^{\ast} (\eta) {\in}
[s_{\nu_i}^{\ast} (\eta) {-} \varepsilon, s_{\nu_i}^{\ast}
(\eta)].$ Since $\varepsilon >0,$ is chosen arbitrarely, we
get: $\lim\limits_{\eta \rightarrow \xi} |
t_{\nu_{i+1}}^{\ast} (\eta) {-}  s_{\nu_i}^{\ast} (\eta)|
{=} 0;$
 but, by the induction hypothesis,
$\lim\limits_{\eta \rightarrow \xi} s_{\nu_{i}}^{\ast}
(\eta)= s_{\nu_{i}}^{\ast} (\xi)= {\tau}_{{i+1}}^{\ast}
(\xi)= t_{\nu_{i+1}}^{\ast} (\xi);$
 hence $t_{\nu_{i+1}}^{\ast} (\cdot)$ is continuous at $\xi,$
as desired.

The fact that  $\eta {\mapsto} z(\eta,t_{\nu_{i+1}}^{\ast}
(\eta))$ is continuous at $\xi$ follows from the estimates
$$ |  z(\eta,t_{\nu_{i+1}}^{\ast} (\eta)) -
z(\xi,t_{\nu_{i+1}}^{\ast} (\xi))| \leq |
z(\eta,t_{\nu_{i+1}}^{\ast} (\eta))-
z(\eta,s_{\nu_{i}}^{\ast} (\eta))| +$$ $$+|
z(\eta,s_{\nu_{i}}^{\ast} (\eta)) -
z(\xi,s_{\nu_{i}}^{\ast} (\xi))| \leq |
\int\limits_{s_{\nu_i}^{\ast}
(\eta)}^{t_{\nu_{i+1}}^{\ast}(\eta)}
\varphi(s,z(\eta,s),v(\eta,s,z(\eta,s)))ds|+$$ $$+ |
z(\eta,s_{\nu_{i}}^{\ast} (\eta)) -
z(\xi,s_{\nu_{i}}^{\ast} (\xi))| \leq \frac{1}{2 L(\eta)}
|t_{\nu_{i+1}}^{\ast} (\eta) - s_{\nu_i}^{\ast} (\eta)|
+$$ $$+| z(\eta,s_{\nu_{i}}^{\ast} (\eta)) -
z(\xi,s_{\nu_{i}}^{\ast} (\xi))|,$$
 from the equalities $ t_{\nu_{i+1}}^{\ast} (\xi) = s_{\nu_i}^{\ast} (\xi) = {\tau}_{i+1}^{\ast} (\xi),$
from the induction hypothesis, and from the continuity of
 $\eta {\mapsto} t_{\nu_{i+1}}^{\ast} (\eta)$ at $\xi.$

To prove that $s_{\nu_{i+1}}^{\ast} (\cdot)$ is continuous
at $\xi,$
 first note, that,
 for some  ${\overline \rho}_1 >0,$ some $\varepsilon>0,$ and
some ${\overline \rho}_2 >0,$ the function ${\tilde t}
(\cdot,\cdot,\cdot, \cdot)$ given by ${\tilde t}
(\eta,t,\tau, y^0) = t - \vartheta_{\nu_{i+1}}
(\eta,y(t,\tau,y^0,v_{r(\nu_{i+1})} ))$ is well defined,
continuous and stictly increasing w.r.t. $t$ for all
$(\eta,t,\tau,y^0)$ in $$ B_{{\overline \rho}_1} (\xi)
\times [s_{\nu_{i+1}}^{\ast} (\xi) {-} \varepsilon, \;
\tau_{i+1}^{\ast} (\xi) {+} \varepsilon] \times
[\tau_{i+1}^{\ast} (\xi) {-} \varepsilon,\;
\tau_{i+1}^{\ast} (\xi) {+} \varepsilon] \times
B_{{\overline \rho}_2} (z(\xi, \tau_{i+1}^{\ast}(\xi)))$$
(the proof is similar to that of  (\ref{r42_65_a})). Then,
since $0= \tilde{t}  \left( \xi, s_{\nu_{i+1}}^{\ast}
(\xi), t_{\nu_{i+1}}^{\ast} (\xi), z(\xi,
t_{\nu_{i+1}}^{\ast}(\xi))\right) ,$ from the implicit
function theorem, it follows that there exists a unique
function
 $(\eta,\tau,y) \mapsto \tilde{s} (\eta,\tau,y)$
that is continuous on the set $${\gothD} {:=} B_{{\tilde
\rho}_1} (\xi) \times  [\tau_{i+1}^{\ast} (\xi) {-}
\tilde{\varepsilon}, \; \tau_{i+1}^{\ast} (\xi) {+}
\tilde{\varepsilon}] \times B_{{\tilde \rho}_2} (z(\xi,
\tau_{i+1}^{\ast}(\xi))), $$ for some ${\tilde \rho}_i \in
]0, {\overline \rho}_i],$ $i{=}1,2,$ $\tilde{\varepsilon}
\in ]0, \varepsilon],$
 and is such that, for every $(\eta,\tau,y^0) {\in} {\gothD}, $
we have: $$   {\tilde s}(\eta,\tau,y^0) =
\vartheta_{\nu_{i+1}} \left(\eta, y \left( {\tilde
s}(\eta,\tau,y^0),\tau,y^0, v_{r(\nu_{i+1})} \right)
\right), \; \; s_{\nu_{i+1}}^{\ast}(\xi) = {\tilde s}
(\xi,t_{\nu_{i+1}}^{\ast}(\xi), z(\xi,
t_{\nu_{i+1}}^{\ast}(\xi))).$$ Then, by the construction,
for each  $\eta$ from some neighborhood of  $\xi {\in}
\rk,$ we get $s_{\nu_{i+1}}^{\ast}(\eta) = {\tilde s}
(\eta,t_{\nu_{i+1}}^{\ast}(\eta), z(\eta,
t_{\nu_{i+1}}^{\ast}(\eta))),$
  and, therefore, having proved that
$\eta \mapsto t_{\nu_{i+1}}^{\ast} (\eta)$ and $\eta
\mapsto z(\eta,t_{\nu_{i+1}}^{\ast} (\eta))$ are continuous
at $\xi,$ we obtain that $ s_{\nu_{i+1}}^{\ast} (\cdot)$ is
continuous at $\xi$ as well.

Finally, the fact that $\eta \mapsto
z(\eta,s_{\nu_{i+1}}^{\ast} (\eta))$ is continuous at
$\xi,$ now follows from the continuity of
$t_{\nu_{i+1}}^{\ast} (\cdot)$ and $\eta \mapsto
z(\eta,t_{\nu_{i+1}}^{\ast}(\eta))$ at $\xi,$ from the
equality $$z(\eta,s_{\nu_{i+1}}^{\ast} (\eta)) =  y \left(
s_{\nu_{i+1}}^{\ast} (\eta), t_{\nu_{i+1}}^{\ast} (\eta),
z(\eta,t_{\nu_{i+1}}^{\ast} (\eta)), v_{r(\nu_{i+1})}
\right),$$ and from theorem on the continuous dependence of
the solution of the Cauchy problem on the initial
condition.

Thus, it is proven by induction over $i {\in}
\{0,1,...,N(\xi) \}$ that (in particular) all
$t_{\nu_i}^{\ast} (\cdot)$ and $s_{\nu_i}^{\ast} (\cdot)$
are continuous at $\xi.$ The proof of lemma 3.4 is
complete.

$\qquad$

Next, we complete the proof of theorem 2.2 as follows.
First, for each $\xi{\in}\rk,$ and for the segment $[t_1,
t_1{+}\sigma(\xi)],$ we get an appropriate approximation of
the original collection of controls $\{ v_{\lambda}
(\xi,\cdot){\}}_{\lambda}$ by another collection of
sufficiently smooth controls $\{ {\hat u}_{\delta_1 ,
\lambda} (\xi,\cdot){\}}_{\lambda}$ that satisfy the
desired boundary conditions (\ref{r43_81}) (see also
condition 5) of theorem 2.2), and, together with their
derivatives w.r.t. $t,$ continuously depend on
$\xi{\in}\rk$ and $t{\in}[t_1, t_1 {+} \sigma(\xi)]$
(lemmas 3.5, 3.6). Second, (lemmas 3.7, 3.8) we do the same
for each $v(\xi,\cdot)$ ($\xi{\in}\rk$) constructed in
lemma 3.4, and for the segment $[t_1{+}\sigma(\xi),T],$ and
obtain a smooth control $v_{\Delta_1} (\xi,\beta,\cdot).$
Finally, for each $\xi{\in}\rk,$ and each
$\beta{\in}\rmpplusl,$ we find the desired control that
satisfies conditions 5), 6) of theorem 2.2 among all
concatenations of $v_{\Delta_1} (\xi,\beta,\cdot)$ with all
possible $ {\hat u}_{\delta_1 , \lambda} (\xi,\cdot) ,$
$\lambda{\in} \overline{B_{\varepsilon_1 (\xi)}(0)} $ (see
(\ref{r43_92})) Then, condition 4) of theorem 2.2 will
follow directly from our construction.

{\bf   Lemma 3.5.}  {\sl  For each
  $ \delta_1(\cdot) {\in} C(\rk; ]0,+ \infty[),$
there exists a family of controls
  $\{  u_{\delta_1}(\xi,\cdot) {\}}_{\xi \in \rk}$
such that,  for each $\xi {\in} \rk,$ we have
$u_{\delta_1}(\xi,\cdot) \in C^1([t_1,t_1 {+}\sigma(\xi)];
\rmpplusl),$ and

1) Functions   $  u_{\delta_1}(\cdot,\cdot) $ and
$\frac{\partial}{\partial t}  u_{\delta_1}(\cdot,\cdot) $
are continuous everywhere on $ \{ (\xi,t) {\in} \rk
{\times}J  | \;   t_1 \leq t \leq t_1 {+}\sigma(\xi)  \}. $

2)  For each $ \xi {\in} \rk,$ we have: $$
u_{\delta_1}(\xi,t_1)= {x}_{p+1}^{\ast}, \; \;
\dot{u}_{\delta_1}(\xi,t_1)=    {z}_{p+1}^{\ast};
\n{r43_81}$$ $$ \max\limits_{t \in [t_1,t_1 +\sigma(\xi)]}
| u_{\delta_1}(\xi,t) - \phi(t,y(\xi,t), \dot{x}_p
(\xi,t))| < \delta_1(\xi).  \n{r43_82}$$ }

If $\sigma(\cdot)$ is a constant function, i.e.,
$t_1{+}\sigma(\xi) {=}t_2$ for some
 $ t_2 \in ]t_1,+ \infty[,$
and for all $\xi{\in}\rk,$ then the proof is the same as
the proof of lemma 3.5 from \cite{ssp_op1}, and follows
from lemma 4.2. If $\sigma(\cdot)$ is not a constant
function, then, taking an arbitrary $t_2{>}t_1,$ and
introducing the linear transformation of time
$t(\xi,s)=t_1+ \frac {\sigma(\xi)(s {-} t_1)}{t_2 {-}
t_1},$ $s{\in}[t_1,t_2],$ $\xi{\in}\rk,$ we reduce the
problem to the case of constant $\sigma(\cdot),$ and
complete the proof of lemma 3.5.

To each  $ \delta_1(\cdot) {\in}C(\rk; ]0,+ \infty[)$
assign a family $\{  u_{\delta_1}(\xi,\cdot) {\}}_{\xi
{\in} \rk},$
 of controls,
 obtained from lemma 3.5, and, then,
for every  $\lambda {=} (\lambda_1,..., \lambda_k)^T  {\in}
\rk,$ and every $\xi {\in} \rk,$ define the control
 $ {\hat u}_{\delta_1,\lambda}(\xi,\cdot) $
on   $[t_1,t_1 {+}\sigma(\xi)] $  by $$ {\hat u}_{\delta_1,
\lambda}(\xi,t) {=} u_{\delta_1}(\xi,t)
{+}\sum\limits_{j=1}^{k} \lambda_j w_j(\xi,t), \; \; \; \;
 \mbox{ whenever } t {\in} [t_1,t_1 {+}\sigma(\xi)]. \n{r43_83} $$
For each  $ \delta_1(\cdot) {\in} C(\rk; ]0,+ \infty[)$
define the family of maps
      $\{  \hat{\Phi}_{\delta_1}(\xi,\cdot) {\}}_{\xi \in \rk}$
 from $\rk$ to $\rk$ as follows:
for each $\xi {\in} \rk,$  and each $\lambda {\in} \rk$
such that $t {\mapsto} y(t, t_1,
y^{\ast},\hat{u}_{\delta_1,\lambda}(\xi, \cdot))$
 is defined for all  $t{\in}[t_1,t_1 {+}\sigma(\xi)], $
by definition, put: $
\hat{\Phi}_{\delta_1}(\xi,\lambda){:=}  y(t_1{+}
\sigma(\xi),t_1, y^{\ast},\hat{u}_{\delta_1,\lambda}(\xi,
\cdot)). $

{\bf   Lemma 3.6.}  {\sl  There exists a function
 $ \delta_1(\cdot) \in C(\rk; ]0,+ \infty[)$   such that
the following conditions hold:

1) For each $(\xi, \lambda) {\in} \Pi_1 {:=} \{(\xi,
\lambda){\in} \rk {\times} \rk |\;  \lambda {\in}
{B_{\varepsilon_1(\xi)}(0)}    \},$
 $t {\mapsto} y(t, t_1, y^{\ast},\hat{u}_{\delta_1,\lambda}(\xi, \cdot))$
is defined for all   $t{\in}[t_1,t_1 {+}\sigma(\xi)], $
and, therefore, $ \hat{\Phi}_{\delta_1}(\xi,\lambda)$ is
well defined.

2) For each  $\xi {\in} \rk,$  the map  $\lambda{\mapsto}
\hat{\Phi}_{\delta_1}(\xi,\lambda)$ is differentiable for
all   $\lambda{\in} {B_{\varepsilon_1(\xi)}(0)} ,   $ and
the maps $( \xi,\lambda) {\mapsto}
\hat{\Phi}_{\delta_1}(\xi,\lambda),$  and
   $( \xi,\lambda)  {\mapsto} \frac{\partial \hat{\Phi}_{\delta_1} }{\partial \lambda}(\xi,\lambda)$
are of classes  $C(\Pi_1; \rk)$   and $C(\Pi_1; \rktk)$
respectively.

3) For each $( \xi, \lambda) {\in} \Pi_1,$ we have: $$ |
\hat{\Phi}_{\delta_1}(\xi,\lambda) -  {\Phi} (\xi,\lambda)
| < \frac{\varepsilon_2(\xi)}{4},    \; \; \; \; \;
\;\mbox{ and } \; \;\;  \;
 \parallel  \frac{\partial \hat{\Phi}_{\delta_1} }{\partial \lambda}   (\xi,\lambda) -  \frac{\partial {\Phi} }{\partial \lambda}  (\xi,\lambda) \parallel  < \rho.  \;  \n{r43_85}$$
} Again, if, for some $t_2{>}t_1,$ we have
$\sigma(\xi){=}t_2$ whenever $\xi{\in}\rk,$ then, the proof
of lemma 3.6 is the same as for lemma 3.6 in
\cite{ssp_op1}. If $\sigma (\cdot)$ is not a constant
function, then, taking some $t_2{>}t_1,$ we easily reduce
the proof to the case of constant $\sigma(\cdot)$ after the
linear transformation of time $t{=}t(\xi,s){=}t_1+ \frac
{\sigma(\xi)(s {-} t_1)}{t_2 {-} t_1},$ $s{\in}[t_1,t_2],$
$\xi{\in}\rk.$

Let $ \delta_1(\cdot) \in C(\rk; ]0,+ \infty[)$ be chosen
from lemma 3.6.

{\bf   Lemma 3.7.}  {\sl For each function $\Delta_1(\cdot,
\cdot) {\in}C(\rk {\times} \rmpplusl; ]0,+ \infty[),$ there
exists a family
        $\{  v_{\Delta_1}(\xi,\beta,\cdot) {\}}_{(\xi,\beta) {\in} \rk {\times} \rmpplusl}$
 of controls of class
 $C^1(J; \rmpplusl)$ such that:

1) The map  $( \xi,\beta)  {\mapsto}
v_{\Delta_1}(\xi,\beta,\cdot)$ is of class  $C(\rk {\times}
\rmpplusl;C^1(J; \rmpplusl)).$

2) For each  $( \xi,\beta) {\in} \rk {\times} \rmpplusl,$
we have: $$v_{\Delta_1}(\xi,\beta,T) {=} \beta;   \; \;
v_{\Delta_1}(\xi,\beta,t_1{+}\sigma(\xi)) {=}
u_{\delta_1}(\xi,t_1{+}\sigma(\xi)),
 \; \;  \dot{v}_{\Delta_1}(\xi,\beta,t_1{+}\sigma(\xi)) {=}   \dot{u}_{\delta_1}(\xi,t_1{+}\sigma(\xi)); \n{r43_86}$$
$${ \parallel v_{\Delta_1}(\xi,\beta,\cdot)  {-} v
(\xi,\cdot)  \parallel }_{L_1(J; \rmpplusl)} <
\Delta_1(\xi,\beta) ;  \n{r43_87}$$ $${ \parallel
v_{\Delta_1}(\xi,\beta,\cdot)   \parallel }_{C(J;
\rmpplusl)} <   2 \max  \left\{|\beta|, \; |
u_{\delta_1}(\xi,t_1{+}\sigma(\xi))| ,\;  M(\xi)  \right\}
{+} 1,  \n{r43_88}$$ where $M(\cdot)$ is defined in lemma
3.4.}

The proof of lemma 3.7 is based on lemma 4.2 and is similar
to the proof of lemma 3.5 from \cite{ssp_op1}.

To each $\Delta_1(\cdot, \cdot) {\in}C(\rk {\times}
\rmpplusl; ]0,+ \infty[),$ we assign a family  $\{
v_{\Delta_1}(\xi,\beta,\cdot) {\}}_{(\xi,\beta) {\in} \rk
{\times} \rmpplusl},$
 of controls,
 obtained from lemma 3.7. Then  the following lemma holds.

{\bf   Lemma 3.8.}  {\sl There exists a function
  $\Delta_1(\cdot, \cdot) {\in}C(\rk {\times} \rmpplusl; ]0,+ \infty[)$
such that, for each $( \xi,\beta) {\in} \rk {\times}
\rmpplusl,$
 the trajectory     $t {\mapsto} y(t, T,\xi, {v}_{\Delta_1}(\xi,\beta, \cdot))$  is defined
for all   $t{\in}J, $  and
 $$ |   y(t, T,\xi, {v}_{\Delta_1}(\xi,\beta, \cdot)) {-} y(t, T, \xi, v(\xi, \cdot))  | {<} \frac {\varepsilon_2 (\xi)}{4}  \; \;  \; \mbox{ whenever }
  t {\in} J, \;   ( \xi,\beta) {\in} \rk {\times} \rmpplusl \n{r43_89}$$
} The proof of lemma 3.8 follows from standard arguments
based on the Gronwall-Bellman lemma.

Let $\Delta_1(\cdot,\cdot)$ be chosen from lemma 3.8. From
lemmas 3.3 and 3.4, we obtain that $ | y(t_1{+}
\sigma(\xi),T, \xi, v(\xi, \cdot)) {-} \Phi(\xi,0)  | {<}
\frac{\varepsilon_2(\xi)}{4} $
 whenever  $\xi {\in} \rk;$
then, using (\ref{r43_89}), we obtain that
 $$ | y(t_1{+} \sigma(\xi),T, \xi, v_{\Delta_1}(\xi,\beta, \cdot)) {-} \Phi(\xi,0)  | {<} \frac{\varepsilon_2(\xi)}{2} \; \; \; \; \; \mbox{ whenever }  (\xi,\beta){\in}\rk{\times}\rmpplusl. \n{r43_90}$$
In addition, from condition $3)$ of lemma 3.6, from
(\ref{r4_23}) (see lemma 3.2),
 from the definition of $\rho,$ and from
 lemma 4.3,
we obtain that, for each    $\xi {\in} \rk, $ the maps $
{\Phi}(\xi,\cdot)$ and
  $  \hat{\Phi}_{\delta_1}(\xi,\cdot)$
are diffeomorphisms of $B_{\varepsilon_1(\xi)}(0)$ onto   $
{\Phi}(\xi,B_{\varepsilon_1(\xi)}(0))$ and onto $
\hat{\Phi}_{\delta_1}(\xi,B_{\varepsilon_1(\xi)}(0))$
respectively. For every    $\xi {\in} \rk, $ by
${\Phi}^{-1}(\xi,\cdot)$ we denote the diffeomorphism (of $
{\Phi}(\xi,B_{\varepsilon_1(\xi)}(0))$  onto   $
B_{\varepsilon_1(\xi)}(0)$)
  that is inverse to
 $ {\Phi}(\xi,\cdot).$ Then, from lemma 3.2,
we obtain that the map
 $\eta {\mapsto}  \hat{\Phi}_{\delta_1} (\xi, {\Phi}^{-1}(\xi,\eta))$
is well defined and continuous at each
 $\eta{\in} \overline{ B_{\varepsilon_2(\xi)}( {\Phi}(\xi,0))}. $
Furthermore, from
  (\ref{r43_85}),
 we get
$ | \eta - \hat{\Phi}_{\delta_1}
(\xi,{\Phi}^{-1}(\xi,\eta))|< \frac{\varepsilon_2(\xi)}{4}$
whenever $ \eta {\in} \overline{ B_{\varepsilon_2(\xi)}(
{\Phi}(\xi,0))}. $ Therefore, using  lemma 4.4, we get
 $ \overline{ B_{\frac{ 3 \varepsilon_2(\xi)}{4}}( {\Phi}(\xi,0))} \subset \hat{\Phi}_{\delta_1}(\xi,B_{\varepsilon_1(\xi)}(0))$
for all   $\xi {\in} \rk. $ From this, and from
(\ref{r43_90}), we obtain that, for each   $( \xi,\beta)
{\in} \rk {\times} \rmpplusl,$ there exists a unique
${\lambda}^\ast (\xi, \beta) {\in}
B_{\varepsilon_1(\xi)}(0)$  such that $$
\hat{\Phi}_{\delta_1}(\xi, {\lambda}^\ast(\xi, \beta)) =
y(t_1 {+} \sigma(\xi), T, \xi, v_{\Delta_1}(\xi,\beta,
\cdot)).  \;\;\;            \n{r43_91}$$ Since family $\{
v_{\Delta_1}(\xi,\beta,\cdot)
{\}}_{(\xi,\beta){\in}\rk{\times}\rmpplusl}$ satisfies
condition 1) of lemma 3.7, the map
  $(\xi,\beta)  {\mapsto} y(t_1{+}\sigma(\xi),T,\xi, v_{\Delta_1}(\xi,\beta,\cdot))$
is of class    $C( \rk {\times} \rmpplusl;\rk);$ hence,
using condition 2) of lemma 3.6 and the implicit function
theorem, we obtain that the map
   $( \xi,\beta)  {\mapsto} {\lambda}^\ast(\xi, \beta)$
is of class  $C( \rk \times \rmpplusl;\rk)$.

  For each  $(\xi,\beta) {\in} \rk {\times} \rmpplusl,$
let $\hat{v}_{(\xi,\beta)}(\cdot)$ be the control given by
$$  \hat{v}_{(\xi,\beta)}(t)  = \left\{ \begin{array}{l}
 \hat{u}_{\delta_1, {\lambda}^\ast(\xi, \beta)}(\xi,t) \; \; \;\; \;\quad\quad  \; \;  \;  \mbox{ if }\;  t_1 \leq t \leq t_1 + \sigma(\xi)  \\
{v}_{\Delta_1}(\xi, \beta,t)     \; \; \; \; \; \; \; \; \;
\; \; \; \quad\quad  \mbox{ if }   \;  t_1 + \sigma(\xi) <
t  \leq  T. \; \;
\end{array}\right.  \n{r43_92} $$

It is clear that the family
     $\{ \hat{v}_{(\xi,\beta)}(\cdot) {\}}_{(\xi,\beta) {\in} \rk {\times} \rmpplusl}$
defined by  (\ref{r43_92}) satisfies conditions 4)-6) of
theorem 2.2. Indeed, condition 6) follows from
(\ref{r43_91}), (\ref{r43_92}),
 and from the definition of
 $  \hat{\Phi}_{\delta_1}(\cdot,\cdot).$
Condition $4)$ follows from  (\ref{r43_86}),
(\ref{r43_83}),
 from condition 1) of lemma
 3.7, from condition 1) of lemma 3.5,
from (\ref{r4_19}), from (\ref{r4_18}), and from the fact
that ${\lambda}^\ast (\cdot, \cdot)$ is continuous.
Finally, condition 5) follows from
 (\ref{r43_81}) (see lemma 3.5),
from (\ref{r43_86}) (see lemma 3.7), and from
(\ref{r4_18}). This completes the proof of theorem 2.2 as
well as the proofs of theorems 2.1, and 1.1-1.3.

\begin{center}
{\bf 4. Appendix.}
\end{center}

{\bf Lemma 4.1.} {\sl Consider a family of control systems
$$   \left\{ \begin{array}{l} \dot z_i (t) =
\sum\limits_{j=1}^{i{+}1} A_{ij} (\xi,t) z_j(t) , \; \; \;
i{=}1,...,p{-}1;\\ \dot z_p (t) = \sum\limits_{j=1}^{p}
A_{pj} (\xi,t) z_j(t) + A_{p \; p{+}1} (\xi,t) w(t),
\end{array} \right.  \; \; \; t\in[t_1, t_1{+}\sigma(\xi)]
\n{add1}   $$ where $\xi{\in}\rbn$ is the parameter of the
family, $\sigma(\cdot)$ is a function of class
$C(\rbn;]0;+\infty[),$ each $A_{ij} (\xi,t)$ is a matrix of
dimension $m_i{\times}m_j$ ($i{=}1,...,p; $
$j{=}1,...,i{+}1$), $w {\in} \rmpplusl$ is the control,
$(z_1,...,z_p{)}^T {\in}\rk {=} \rmlplusmp$ is the state,
and $z_i {\in} \rmi,$ $i{=}1,...,p.$

Assume that all $A_{ij} (\cdot,\cdot)$ are continuous  on
$D:= \{ (\xi,t) {\in} \rbn {\times} \rr | \;
t_1{\leq}t{\leq}t_1{+}\sigma(\xi)  \}.$ Suppose that, for
each $i{=}1,...,p,$ we have $m_i{\leq}m_{i{+}1},$ and there
exist numbers $j_1{<}j_2{<}...{<}j_{m_i}$ in $\{
1,...,m_{i{+}1}  \}$ such that the columns $a_{i\;
i{+}1}^{j_1} (\cdot,\cdot),...,a_{i\; i{+}1}^{j_{m_i}}
(\cdot,\cdot),$ of matrix $A_{i\; i{+}1} (\cdot,\cdot),$
numbered by $j_1,...,j_{m_i}$ satisfy the condition ${\rm
det} [a_{i \; i{+}1}^{j_1} (\xi,t),..., a_{i\;
i{+}1}^{j_{m_i}} (\xi,t)] \not=0$ for all $(\xi,t){\in}D.$

Then, for each $z^T {\in} \rk,$ and each $\mu{\in}\nn,$
there exists a family of controls $\{ w(\xi,\cdot)
{\}}_{\xi{\in}\rbn}$ such that

1) For each $(\xi,t){\in}D,$ and each $l{=}0,1...,{\mu},$
there exist $\frac{{\partial}^l w}{\partial {t}^l} (\xi,t)$
and $\frac{{\partial}^l w}{\partial {t}^l} (\cdot,\cdot)
{\in} C(D; \rmpplusl),$ $(l{=}0,1,...,\mu).$

2) For each $\xi {\in} \rbn,$ the control $w(\xi,\cdot)$
steers $0{\in}\rk$ into $z^T$ in time $[t_1, t_1{+}
\sigma(\xi)]$ w.r.t. (\ref{add1}).

3) For each $\xi{\in}\rbn,$ we have: $\frac{{\partial}^l
w}{\partial t^l} (\xi,t_1){=} \frac{{\partial}^l
w}{\partial t^l} (\xi,t_1{+}\sigma(\xi)){=}0{\in}
\rmpplusl,$ $l{=}0,1,...,\mu.$ }

If $\sigma(\cdot)$ is a constant function, i.e.,
$[t_1,t_1{+}\sigma(\xi)] = [t_1,t_2]$ for some  $t_2{>}t_1$
and for all $\xi \in \rbn,$ then, the proof of lemma 4.1 is
similar to the proof of lemma 3.1 in \cite{ssp_op1} and is
based on the same construction as in theorem 1 (sect. 2) of
\cite{pavl3}.
If $\sigma(\cdot)$ is  not a constant function, then,
taking some arbitrary $t_2{>}t_1,$ and introducing the
linear transformation of time $t(\xi,s) = t_1 +
\frac{\sigma(\xi) (s-t_1)}{t_2-t_1},$ $s{\in}[t_1,t_2],$
$\xi{\in}\rbn,$ we reduce the problem to the case of
constant $\sigma(\cdot).$

{ \bf Lemma 4.2.} {\sl Let  $X$ and $Y$ be linear normed
spaces equipped with norms
 $\parallel \cdot \parallel_X$ and
 $ \parallel \cdot \parallel_Y$ respectively.
 Assume that  $Y \subset X, $ and
 $Y$ is dense in $X$  w.r.t. $\parallel \cdot \parallel_X.$ Suppose that
 $\neta \mapsto f(\neta)$
  is of class  $C(\rbn;X)$, where  $N~\in~\nn.$
Then:

A) For each  $\nrho >0$  there exists a map   $\neta
\mapsto \hat{f}(\neta)$
 of class $C(\rbn;Y)$   such that  for every  $\neta \in \rbn$ we have:
$$ {  \parallel  f(\neta)  -  \hat{f}(\neta)\parallel}_X  <
\nrho.   \n{r13_1}$$

B) For each function $\neta \mapsto \nrho( \neta) $   of
class $C(\rbn;]0,+\infty[)$
   there exists a map  $\neta \mapsto \hat{g}(\neta)$    of  class   $C(\rbn;Y)$
   such that for every  $\neta \in \rbn$   we have:
$$ {  \parallel  f(\neta)  -  \hat{g}(\neta)\parallel}_X  <
\nrho( \neta).   \n{r13_2}$$ } {\bf Proof of lemma 4.2.}
A).
    To each $\neta \in \rbn$  assign an open ball
  $B_\neta$  with the center at $\neta \in \rbn$ such that for each
$\overline{\neta} \in B_\neta$
    we have:
$ {  \parallel  f(\overline{\neta})  -
f(\neta)\parallel}_X  < \frac{\nrho}{2}.   $
   Let  $ \{\varphi_i(\cdot)\}_{i=1}^\infty $
      be the partition
     of unity corresponding to the family
  of open sets $\{ B_{\neta} {\}}_{\neta \in \rbn},$
     which covers $\rbn;$  i.e.,
     ~\cite[p.~66]{sternberg},
  $ \{\varphi_i(\cdot)\}_{i=1}^\infty $ is a sequence of functions, of class
  $C^{\infty}(\rbn; \rr),$ having compact supports
   such that the following conditions hold:

    1)   $\{ {\rm supp} \;  \varphi_i (\cdot)  {\}}_{i=1}^{\infty}$
      is a  locally finite covering of  $\rbn$
    (in the sense that each compact set in $\rbn$
     has the nonempty intersection with only a finite number of the supports),
   and for each $i \in \nn$
    and each $\neta \in \rbn$    we have: $0 \leq  \varphi_i(\neta) \leq 1;$

    2) For each $i \in \nn$ there exists $\neta_{i} \in \rbn$
     such that ${\rm supp} \; \varphi_i(\cdot) \subset B_{\neta_{i}};$

    3) For each $\neta \in \rbn$  we have: $ \sum_{i=1}^{\infty}\varphi_i(\neta)  =1.$

    (In condition 3), for each fixed $\neta $ the sum is finite due to 1)).
     For each $i \in \nn$  fix   $y_i \in Y$  such that
    $ {  \parallel  y_i  -  f(\neta_{i})\parallel}_X  < \frac{\nrho}{2}.   $
     Then, the map
$\neta \mapsto \hat{f}(\neta) =
\sum\limits_{i=1}^{\infty}\varphi_i(\neta)y_i    $
    satisfies (\ref{r13_1}). Indeed, using conditions 1)-3),
    we obtain:
$$   \parallel  \hat{f}(\neta)  -  f(\neta) \parallel_X  =
\parallel \sum\limits_{i=1}^{\infty}\varphi_i(\neta)y_i - (
\sum\limits_{i=1}^{\infty}\varphi_i(\neta))f(\neta)
\parallel_X = \parallel
\sum\limits_{i=1}^{\infty}\varphi_i(\neta)(y_i   - f(\neta)
)  \parallel_X \leq $$
 $$ \leq   \sum\limits_{i=1}^{\infty}\varphi_i(\neta) \parallel y_i   - f(\neta_{i})   {\parallel}_X
+  \sum\limits_{i=1}^{\infty}\varphi_i(\neta){\parallel
f(\neta_{i})   - f(\neta)   \parallel}_X  \leq
\frac{\nrho}{2}+ \frac{\nrho}{2}  = \nrho.$$
    The fact that the map  $\neta \mapsto \hat{f}(\neta)$
  is of class   $C(\rbn; Y)$  follows from
  the inclusion   $\varphi_i(\cdot) \in C^{\infty}(\rbn; \rr)$
   and from 1). A) is now proven.

    Let us prove B).  For each  $k \in \nn,$ put: $\nrho_{k} = \min \limits_{\parallel \neta\parallel \leq k}{ \nrho(\neta)}$.
    Then, from A), it follows that for each
     $k \in \nn$ there exists a map   $\neta \mapsto \hat{f_k}(\neta)$
      of class $C(\rbn; Y)$  such that for each $\neta \in \rbn$   we have:
   $   \parallel  {\hat{f}}_k(\neta)  -  f(\neta) {\parallel}_X \leq \nrho_{k}.$
    Let    $\neta~\mapsto~\hat{g}(\neta)$   be the map of $\rbn$  to $Y$
    given by
   $$\hat{g}(\neta) = (k- | \neta|) {\hat f}_k (\neta)  + (1-k+ | \neta|) {\hat f}_{k+1} (\neta),
   \mbox { whenever } k-1 \leq  | \neta| < k, \; \; k\in \nn. $$
    Then  $\neta \mapsto \hat{g}(\neta)$   is of class  $C(\rbn;Y)$
      and satisfies  (\ref{r13_2}). The proof of lemma 4.2 is complete.

    {\bf Lemma 4.3.} {\sl Assume that  $B  \subset \rk$
    is a convex open set, and
     $F(\cdot) \in C^{1}(B; \rk)$  satisfies the condition:
     for each  $\lambda_0 \in B $ the matrix
   $\frac{\partial F}{\partial \lambda}(\lambda_0)$  is positive definite.
    Then,
   $\lambda \mapsto F(\lambda)$ is a diffeomorphism of
     $B$  onto $F(B)$.}

    { \bf Proof of lemma 4.3.}
     Take any  $\lambda_1 \in B$ and any $\lambda_2 \in B$ such that
$\lambda_1 \not= \lambda_2 .$
    It is sufficient to prove that
  $\langle F(\lambda_1) -  F(\lambda_2), \; \;  \lambda_1 -  \lambda_2 \rangle  \not=0,$
  (where $\langle \cdot,\cdot \rangle$  is the scalar product in $\rk$) which
   implies   $ F(\lambda_1) \not= F(\lambda_2).$
   Then, from the implicit function theorem, we will obtain that
     $F^{-1}$ is of class $  C^{1}(F(B); \rk).$
     From the assumptions it follows that
$$\langle F(\lambda_1) -  F(\lambda_2), \; \;  \lambda_1 -
\lambda_2 \rangle = \langle  F(\lambda_2 + \theta
(\lambda_1 -  \lambda_2 ))
 -  F(\lambda_2), \; \;  \lambda_1 -  \lambda_2  \rangle |_{\theta=0}^{\theta=1} = $$
$$= \int\limits_{0}^{1}\langle  \frac {\partial F}
{\partial\lambda} (\lambda_2 + \theta(\lambda_1 -
\lambda_2  )) (\lambda_1 -  \lambda_2),  \;\; \lambda_1 -
\lambda_2 \rangle d\theta > 0,$$
 which completes the proof of the lemma.

{\bf Lemma 4.4.} {\sl Let $f(\cdot)$ be a continuous map of
       $B=\{ x \in \rn| \; |x-x^0| \leq r    \}$  to $\rn,$  where
$r>0.$
    Assume that there exists
$\varepsilon \in ]0,r[$ such that
   $|f(x) - x| \leq r-\varepsilon$ for all
   $x \in B.$ Then, each point
    $z \in \rn $  such that $|z-x^0| \leq \varepsilon $
    belongs to the image of $B$ under
   $f,$  i.e., there exists
    $x^{\ast} \in B$ such that  $f(x^{\ast}) =z.$
}

    The lemma follows from the Brouwer fixed point theorem.
    The proof is given in ~\cite[p.~276-~277]{lee_markus}.



\end{large}
\end{document}